\documentclass[11pt]{article}
\usepackage[T1]{fontenc}
\usepackage{lmodern,amsmath,amssymb,amsthm,mathtools,booktabs,array}
\usepackage[a4paper,margin=25mm]{geometry}
\usepackage{xcolor}

\usepackage{tikz}
\usetikzlibrary{backgrounds,fit}
\usepackage{placeins,needspace}
\usepackage[hidelinks]{hyperref}
\newtheorem{problem}{Problem}
\newcommand{\Z}{\mathbb Z}
\newtheorem{theorem}{Theorem}[section]
\newtheorem{lemma}[theorem]{Lemma}
\newtheorem{proposition}[theorem]{Proposition}
\newtheorem{corollary}[theorem]{Corollary}
\theoremstyle{definition}
\newtheorem{definition}[theorem]{Definition}

\title{Connected regular $\Gamma$-distance magic graphs with magic constant zero}
\author{{ Ahmet Batal$^{1,2}$\thanks{Ahmet Batal's research was partially supported by program ``Excellence initiative -- research university'' (IDUB) for the AGH University.}, {Sylwia Cichacz$^{2}$\thanks{Sylwia Cichacz's research was supported by  the AGH University of Krakow under grant no. 16.16.420.054, funded by the Polish Ministry of Science and Higher Education.}, Hubert Grochowski$^3$}}
	\footnote{e-mails: cichacz@agh.edu.pl, ahmetbatal@iyte.edu.tr, hubert.grochowski@pw.edu.pl}\\
\normalsize $^1$Department of Mathematics, Izmir Institute of Technology, 35430, Urla, Izmir, Turkey\\
\normalsize $^2$AGH University of Krakow, Poland \\
\normalsize $^3$Warsaw University of Technology, Poland\\
}
\date{\small 16 October 2026}
\begin{document}
\maketitle
\begin{abstract}
A graph of order $n$ is distance magic if it admits a bijective
vertex labeling by $\{1,\ldots,n\}$ such that the sum of the labels
of the neighbours is the same at every vertex.
Let $(\Gamma,+)$ be a finite Abelian group. A $\Gamma$-distance magic
labeling of a graph is a bijection from its vertex set to $\Gamma$
such that the sum of the labels of the neighbours of each vertex
is a fixed element of $\Gamma$, called the magic constant.
The constructions of Kov\'a\v{r} and collaborators provide regular
distance magic graphs of odd order for every even degree from $4$ to $14$.
We adapt these constructions to prove the existence of connected
$r$-regular $(A\times\mathbb{Z}_n)$-distance magic graphs with magic
constant $0$, where $A$ is any finite Abelian group of odd order,
$n$ is odd, and $r\in\{4,6,8,10,12,14\}$, under explicit sufficient
conditions on $n$. Here $\mathbb{Z}_n$ is the cyclic group of order $n$.
It suffices to take $n\ge13$ for $r=6$, $n\ge19$ for $r=14$,
and $n\ge15$ for $r\in\{4,8,10,12\}$.
\end{abstract}

\noindent\textbf{2020 Mathematics Subject Classification.}
Primary 05C78; Secondary 05C25, 05C07.

\smallskip
\noindent\textbf{Keywords.}
Group distance magic labeling, regular graph, finite Abelian group,
magic constant zero, recursive construction.

\section{Introduction}
All graphs considered are finite, simple and undirected.
Let $G$ be a graph of order $n$. For a vertex $v\in V(G)$, we write
$N_G(v)$ for its set of neighbours and $d_G(v)=|N_G(v)|$ for its degree.
For a vertex labeling $f$ of $G$ with integer values or values in an
Abelian group, the \emph{weight} of $v$ is
\[
 w_f(v)=\sum_{u\in N_G(v)}f(u),
\]
where the sum is taken in the corresponding additive group.
A \emph{distance magic labeling} of $G$ is a
bijection $f\colon V(G)\to\{1,\ldots,n\}$ for which all vertex
weights are equal. If the graph is $r$-regular, summing the weights
over all vertices shows that their common value is $r(n+1)/2$.
For background, see~\cite{ArumugamSurvey,Gallian}.

The existence problem for regular distance magic graphs asks which
orders and degrees admit such a labeling. The case of positive
degree and even order was settled in~\cite{FroncekFair}.
For odd order, the degrees $4$, $6,8,10,12$, and $14$ were treated
in~\cite{KFK4reg,KSKK6810,KovKrb14}, respectively.

Replacing the consecutive integer labels by the elements of a finite
Abelian group gives the notion of group distance magic
labeling~\cite{FroncekGroup}. Let $(\Gamma,+)$ be such a group, written
additively, and suppose that $|V(G)|=|\Gamma|$.
A bijection $f\colon V(G)\to\Gamma$ is a
\emph{$\Gamma$-distance magic labeling} with \emph{magic constant}
$\mu\in\Gamma$ if $w_f(v)=\mu$ for every vertex $v$.
We call $G$ \emph{$\Gamma$-distance magic} if it admits such a labeling.

We study connected regular graphs with a group distance magic labeling
whose magic constant is $0$. Our label groups have the form
$A\times\Z_n$, where $A$ is a finite Abelian group of odd order and
$\Z_n=\Z/n\Z$ is the cyclic group of odd order $n$.
An integer distance magic labeling of an $r$-regular graph of odd
order $n$ yields a $\Z_n$-distance magic labeling with constant $0$
by subtracting $(n+1)/2$ from every label and reducing modulo $n$.
The question addressed here is how to enlarge such a graph so that
the labels exhaust $A\times\Z_n$, while keeping the degree fixed,
the graph connected, and all weights zero.

The recursive constructions of~\cite{KFK4reg,KSKK6810,KovKrb14}
provide the starting point. We adapt their edge exchanges to group
labels and add constructions that complete the labeling for every
odd order of $A$. The resulting existence theorem is as follows.

\begin{theorem}\label{thm:main}
Let $A$ be a finite Abelian group of odd order, and let $n$ be
a positive odd integer.

For every $r\in\{4,8,10,12\}$, there exists a connected
$r$-regular $(A\times\Z_n)$-distance magic graph with magic
constant $0$ whenever $n\ge15$.
The same conclusion holds for $r=14$ whenever $n\ge19$.

For $r=6$, such a graph exists whenever $n\ge13$.
It also exists when $n=9$ and $|A|\equiv1\pmod4$, or when
$n=11$ and $|A|\equiv3\pmod4$.
\end{theorem}

The theorem does not assert that the excluded values are impossible.
The theorem allows an arbitrary finite Abelian group $A$ of odd
order, including the trivial group.
Its proof begins with a graph labeled by $\{0\}\times\Z_n$ and
adds vertices labeled by the remaining group elements in pairs
$a,-a$. Such a pair contributes zero to a neighbourhood sum whenever
both of its vertices are adjacent to the vertex in question.
The distinction between $|A|\equiv1\pmod4$ and $|A|\equiv3\pmod4$
comes from the number of pairs available for this construction.
The required modifications and final edge exchanges are developed
after the recursive steps.

\section{Preliminaries and base graphs}
\begin{definition}
Let $G$ be a graph of order $n$. A bijection
$g\colon V(G)\to\{1,\ldots,n\}$ is \emph{centered} if
\[
 w_g(v)=d_G(v)\frac{n+1}{2}
\]
for all $v\in V(G)$.
\end{definition}
An integer distance magic labeling of a regular graph is centered.
The definition also applies to the nonregular base graphs used below.

\begin{lemma}\label{lem:reduction}
If $g$ is a centered labeling of a graph of odd order $n$, then
\[
 \ell(v)=g(v)-\frac{n+1}{2}\pmod n
\]
is a $\Z_n$-distance magic labeling with magic constant $0$.
Moreover, $g(u)+g(v)=n+1$ implies $\ell(u)+\ell(v)=0$.
\end{lemma}
\begin{proof}
The shifted labels form a complete residue system modulo $n$, and
\[
 \sum_{u\in N_G(v)}\ell(u)
 =w_g(v)-d_G(v)\frac{n+1}{2}=0\pmod n.
\]
The assertion about pairs follows by subtracting $n+1$ from their sum.
\end{proof}

For a finite Abelian group $A$ of odd order and a positive odd integer
$n$, put $\Gamma=A\times\Z_n$, $K=\{0\}\times\Z_n$, and
$\Omega=\Gamma\setminus K$. Composing the labeling in
Lemma~\ref{lem:reduction} with $x\mapsto(0,x)$ gives an injective
labeling in $\Gamma$ whose image is $K$ and whose weights are zero.
The labels outside $K$ can be added in pairs by the following observation.

\begin{lemma}\label{lem:pairs}
The set $\Omega$ is a disjoint union of
$n(|A|-1)/2$ two-element sets of the form $\{a,-a\}$.
\end{lemma}
\begin{proof}
Since $K$ is a subgroup, $\Omega$ is closed under inversion.
The group $\Gamma$ has odd order, so its only element equal to its
inverse is $0$, which belongs to $K$. The number of pairs is
$|\Omega|/2=n(|A|-1)/2$.
\end{proof}

In the constructions below, the letters $U,V,W,X,Y,Z$ and
$P,Q,T,D$ denote specified two-element sets of vertices. Distinct such pairs
in a construction are disjoint. We write $UV$ for the set of all
four edges with one endpoint in $U$ and the other in $V$, and call
this set a \emph{join}. Under a group labeling, a vertex pair is
\emph{complementary} if its labels sum to zero.

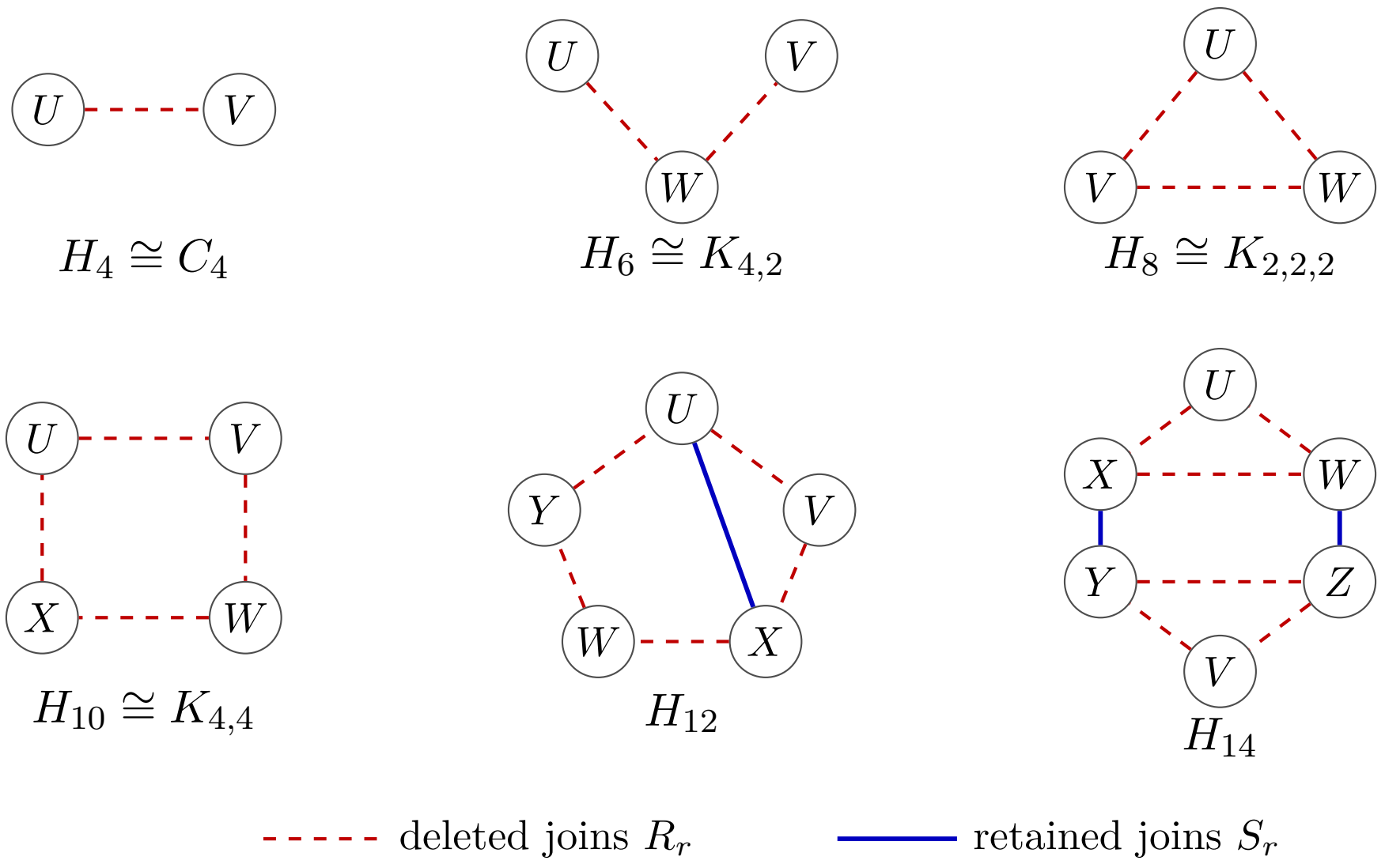
\begin{figure}[!htbp]
\centering
\begin{tikzpicture}[
  x=1cm,y=1cm,
  pair/.style={circle,draw=black!70,fill=white,minimum size=6mm,
    inner sep=1pt,font=\small},
  removed/.style={red!75!black,dashed,line width=.8pt},
  retained/.style={blue!75!black,line width=1.1pt}]
\begin{scope}[shift={(0,0)}]
  \node[pair] (U) at (-.8,0) {$U$};
  \node[pair] (V) at (.8,0) {$V$};
  \draw[removed] (U)--(V);
  \node at (0,-1.25) {$H_4\cong C_4$};
\end{scope}
\begin{scope}[shift={(4.5,0)}]
  \node[pair] (U) at (-1,.45) {$U$};
  \node[pair] (V) at (1,.45) {$V$};
  \node[pair] (W) at (0,-.65) {$W$};
  \draw[removed] (U)--(W)--(V);
  \node at (0,-1.25) {$H_6\cong K_{4,2}$};
\end{scope}
\begin{scope}[shift={(9,0)}]
  \node[pair] (U) at (0,.55) {$U$};
  \node[pair] (V) at (-1,-.65) {$V$};
  \node[pair] (W) at (1,-.65) {$W$};
  \draw[removed] (U)--(V)--(W)--(U);
  \node at (0,-1.25) {$H_8\cong K_{2,2,2}$};
\end{scope}
\begin{scope}[shift={(0,-3.4)}]
  \node[pair] (U) at (-.85,.65) {$U$};
  \node[pair] (V) at (.85,.65) {$V$};
  \node[pair] (W) at (.85,-.85) {$W$};
  \node[pair] (X) at (-.85,-.85) {$X$};
  \draw[removed] (U)--(V)--(W)--(X)--(U);
  \node at (0,-1.65) {$H_{10}\cong K_{4,4}$};
\end{scope}
\begin{scope}[shift={(4.5,-3.4)}]
  \node[pair] (U) at (0,.9) {$U$};
  \node[pair] (V) at (1.15,.05) {$V$};
  \node[pair] (X) at (.7,-1.05) {$X$};
  \node[pair] (W) at (-.7,-1.05) {$W$};
  \node[pair] (Y) at (-1.15,.05) {$Y$};
  \draw[removed] (U)--(V)--(X)--(W)--(Y)--(U);
  \draw[retained] (U)--(X);
  \node at (0,-1.65) {$H_{12}$};
\end{scope}
\begin{scope}[shift={(9,-3.4)}]
  \node[pair] (U) at (0,1.1) {$U$};
  \node[pair] (X) at (-1,.35) {$X$};
  \node[pair] (W) at (1,.35) {$W$};
  \node[pair] (Y) at (-1,-.55) {$Y$};
  \node[pair] (Z) at (1,-.55) {$Z$};
  \node[pair] (V) at (0,-1.3) {$V$};
  \draw[removed] (U)--(X)--(W)--(U);
  \draw[removed] (V)--(Y)--(Z)--(V);
  \draw[retained] (X)--(Y) (W)--(Z);
  \node at (0,-1.85) {$H_{14}$};
\end{scope}
\draw[removed] (1,-6.1)--(2,-6.1);
\node[right,font=\small] at (2,-6.1) {deleted joins $R_r$};
\draw[retained] (5.8,-6.1)--(6.8,-6.1);
\node[right,font=\small] at (6.8,-6.1) {retained joins $S_r$};
\end{tikzpicture}
\caption{Schematic representations of all six subgraphs $H_r$.
Each node represents a prescribed pair of nonadjacent vertices in
$H_r$, and each line represents the four edges of a complete join
$K_{2,2}$. Dashed red joins form $R_r$, and solid blue joins form $S_r$.
The selected copy need not be induced in the ambient graph.}
\label{fig:Hr-schematic}
\end{figure}

For $r\in\{4,6,8,10,12,14\}$, let $H_r$ be the graph with the
vertex pairs listed in Table~\ref{tab:patterns}. Let $R_r$ and $S_r$
be the edge sets consisting of the joins in the corresponding
columns, and set $E(H_r)=R_r\cup S_r$. These edge sets are disjoint.
Figure~\ref{fig:Hr-schematic} gives the corresponding diagrams on
the vertex pairs. In particular, $H_4\cong C_4$, $H_6\cong K_{4,2}$,
$H_8\cong K_{2,2,2}$, and $H_{10}\cong K_{4,4}$.
The graph $H_{12}$ is obtained from a $5$-cycle with a chord by
replacing each vertex with an independent pair and each edge with
a join. For $H_{14}$, the schematic graph on the six pairs consists
of two triangles linked by two edges.

\begin{table}[htbp]
\centering
\caption{The graphs $H_r$. A recursive step deletes the joins in $R_r$
and retains those in $S_r$.}\label{tab:patterns}
\begin{tabular}{clll}
\toprule
$r$ & Pairs of $H_r$ & $R_r$ & $S_r$\\
\midrule
4 & $U,V$ & $UV$ & $\varnothing$\\
6 & $U,V,W$ & $UW,VW$ & $\varnothing$\\
8 & $U,V,W$ & $UV,VW,WU$ & $\varnothing$\\
10 & $U,V,W,X$ & $UV,VW,WX,XU$ & $\varnothing$\\
12 & $U,V,W,X,Y$ & $UV,VX,XW,WY,YU$ & $UX$\\
14 & $U,V,W,X,Y,Z$ & $UX,UW,XW,VY,VZ,YZ$ & $XY,WZ$\\
\bottomrule
\end{tabular}
\end{table}

\begin{definition}\label{def:special}
Let $G$ be a graph with a labeling $f\colon V(G)\to\Lambda$ in an
Abelian group $\Lambda$. A \emph{special copy of $H_r$} in $(G,f)$
is a subgraph isomorphic to $H_r$, together with the prescribed
vertex pairs in Table~\ref{tab:patterns}, such that every prescribed
pair is complementary.
\end{definition}
The copy need not be induced. In particular, vertices belonging to
one prescribed pair may be adjacent in $G$, although they are not
adjacent in the selected subgraph $H_r$.

\begin{proposition}\label{prop:bases}
Let $n$ be odd and let $r\in\{4,6,8,10,12,14\}$.
Under the conditions in the following table, there is a connected
$r$-regular graph $G$ of order $n$ with a distance magic labeling
$g$ and a specified copy of $H_r$ whose prescribed pairs have label
sum $n+1$.
\[
\begin{array}{c|c}
r&n\\\hline
4&n\ge17\\
6&n=9\text{ or }n\ge13\\
8,10,12&n\ge15\\
14&n\ge19.
\end{array}
\]
\end{proposition}
\begin{proof}
For $r=4$, use the initial graph in the proof of Theorem~2.2
and the recursive construction of Lemma~2.1 of~\cite{KFK4reg}.
For $r=6,8,10,12$, use the initial graphs in the proofs of
Theorems~1--4 and the corresponding Lemmas~1--4 of~\cite{KSKK6810}.
For $r=14$, use the initial graphs in the proof of Theorem~3.1
and Lemma~2.2 of~\cite{KovKrb14}.
In each case, the initial graphs contain the specified subgraph with
pair sums equal to the order plus one, and the recursive construction
preserves this property. The pairs may be renamed to agree with
Table~\ref{tab:patterns}.
\end{proof}

\begin{figure}[!htbp]
\centering
\begin{tikzpicture}[
  scale=.75,
  font=\small,
  vert/.style={circle,fill=black,inner sep=0pt,minimum size=3.5pt},
  gedge/.style={draw=black!65,line width=.4pt},
  redge/.style={draw=red!75!black,line width=1pt}
]
\node[vert] (v3)  at ( 1.05, 1.05) {};
\node[vert] (v6)  at (-1.05, 1.05) {};
\node[vert] (v12) at ( 1.05,-1.05) {};
\node[vert] (v9)  at (-1.05,-1.05) {};

\foreach \v/\k in {1/0,14/1,4/2,10/3,11/4,7/5,
                    2/6,5/7,8/8,13/9,0/10}{
  \node[vert] (v\v) at ({90+360*\k/11}:3.25) {};
}

\foreach \u/\v in {
  0/1,0/4,0/8,
  1/6,1/10,1/14,
  2/5,2/7,
  4/9,4/10,4/11,
  5/7,5/8,5/13,
  6/14,7/11,7/12,8/12,8/13,
  9/11,10/11,10/14,13/14}{
  \draw[gedge] (v\u)--(v\v);
}
\draw[gedge] (v0) to[bend left=20] (v2);
\draw[gedge] (v2) to[bend right=14] (v3);
\draw[gedge] (v3) to[bend left=18] (v13);

\draw[redge] (v3)--(v6)--(v12)--(v9)--(v3);

\foreach \v/\k in {1/0,14/1,4/2,10/3,11/4,7/5,
                    2/6,5/7,8/8,13/9,0/10}{
  \node at ({90+360*\k/11}:3.49) {$\v$};
}
\node[above right,fill=white,inner sep=1pt] at (v3)  {$3$};
\node[above left]  at (v6)  {$6$};
\node[below right,fill=white,inner sep=1pt] at (v12) {$12$};
\node[below left]  at (v9)  {$9$};
\end{tikzpicture}
\caption{A $4$-regular $\mathbb Z_{15}$-distance magic graph with magic constant $0$.
The special cycle $3,6,12,9$ is highlighted in red.
Its opposite labels satisfy $3+12=6+9=0$ in $\mathbb Z_{15}$.}
\label{fig:Z15-distance-magic}
\end{figure}
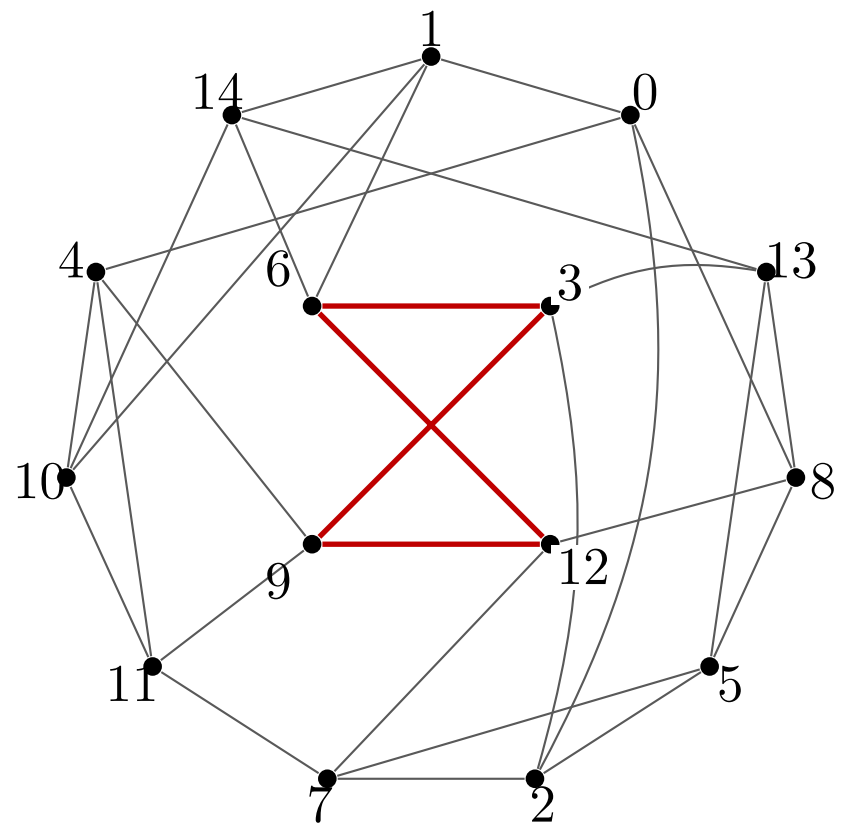

For $r=4$ we also use the connected graph in Figure~\ref{fig:Z15-distance-magic}.
Its vertices are labeled by the distinct residues $0,\ldots,14$.
Each vertex has  degree $4$ and a weight equal to $0$ in $\Z_{15}$.
The highlighted cycle $3,6,12,9$ is the special $H_4$. Its prescribed opposite
pairs are $\{3,12\}$ and $\{6,9\}$, both of sum $0$ in $\Z_{15}$.

\FloatBarrier
\Needspace{16\baselineskip}
\begin{corollary}\label{lem:special-base}
Let $n$ be a positive odd integer and let $r\in\{4,6,8,10,12,14\}$.
There exists a connected $r$-regular graph of order $n$ admitting a
$\Z_n$-distance magic labeling with magic constant $0$ and containing
a special copy of $H_r$ whenever
\[
\begin{array}{c|c}
r&n\\\hline
4,8,10,12&n\ge15\\
6&n=9\text{ or }n\ge13\\

14&n\ge19.
\end{array}
\]
\end{corollary}
\begin{proof}
For $(r,n)\ne(4,15)$, take the graph and integer labeling in
Proposition~\ref{prop:bases}. The labeling is centered because the
graph is regular. Lemma~\ref{lem:reduction} therefore gives a
$\Z_n$-distance magic labeling with magic constant $0$ and makes
every prescribed pair of the specified $H_r$ complementary.
For $(r,n)=(4,15)$, use the graph in
Figure~\ref{fig:Z15-distance-magic}.
\end{proof}

\section{Recursive steps}
Throughout this section, $\Lambda$ is a finite Abelian group.
The recursive steps enlarge a labeled graph by adding vertex pairs
and modifying joins. Unused elements of $\Lambda$ provide the new
labels, so intermediate labelings need only be injective.

Each added or deleted join changes the degree of every incident
vertex by two. For the weight calculations, we use the fact that
a complementary pair contributes zero whenever both of its vertices
belong to a neighbourhood. The following lemma states the conditions
under which these changes preserve all vertex weights.

\begin{lemma}\label{lem:replacement}
Let $f\colon V(G)\to\Lambda$ be injective, with all vertex weights
zero. Modify the edges of $G$ and, if required, add new vertex
pairs labeled by pairwise disjoint two-element inverse pairs outside
$f(V(G))$. Suppose that, at each old vertex, the sets of deleted
and added neighbours are both disjoint unions of complementary pairs.
Suppose also that the neighbourhood of each new vertex is a
disjoint union of complementary pairs. Then the extended labeling is
injective and all its vertex weights are zero.
\end{lemma}
\begin{proof}
Every deleted or added pair contributes zero to the weight of an
old vertex. The same holds for every pair in the neighbourhood of
a new vertex. The choice of the new labels preserves injectivity.
\end{proof}

\begin{figure}[p]
\centering
\begin{tikzpicture}[
  x=1.18cm,y=1.08cm,
  oldpair/.style={circle,draw=black!70,fill=white,minimum size=4.5mm,
    inner sep=.5pt,font=\footnotesize,text=black},
  newpair/.style={oldpair,draw=green!45!black,fill=green!10},
  removed/.style={red!75!black,dashed,line width=.65pt},
  retained/.style={blue!75!black,line width=1.2pt},
  added/.style={green!45!black,line width=.65pt}]
\node[font=\small\bfseries] at (0,1.7) {Before};
\node[font=\small\bfseries] at (5.6,1.7) {After};

\begin{scope}[shift={(0,-0.0)}]
\node[anchor=east] at (-1.6,0) {$r=4$};
\node[oldpair] (bU) at (-0.75,0) {$U$};
\node[oldpair] (bV) at (0.75,0) {$V$};
\node[oldpair] (aU) at (4.3,0.22) {$U$};
\node[oldpair] (aV) at (4.3,-0.22) {$V$};
\node[newpair] (aP) at (6.8,0) {$P$};
\begin{scope}[on background layer]
\draw[removed] (bU) -- (bV);
\draw[added] (aU) -- (aP);
\draw[added] (aV) -- (aP);
\end{scope}
\draw[->,thick] (1.4,0)--(2.6,0);
\end{scope}
\begin{scope}[shift={(0,-2.7)}]
\node[anchor=east] at (-1.6,0) {$r=6$};
\node[oldpair] (bU) at (-0.8,0.55) {$U$};
\node[oldpair] (bV) at (0.8,0.55) {$V$};
\node[oldpair] (bW) at (0,-0.65) {$W$};
\node[oldpair] (aU) at (4.3,0.44) {$U$};
\node[oldpair] (aV) at (4.3,-0.44) {$V$};
\node[oldpair] (aW) at (4.3,0.0) {$W$};
\node[newpair] (aP) at (6.8,.8) {$P$};
\node[newpair] (aQ) at (6.8,-.8) {$Q$};
\begin{scope}[on background layer]
\draw[removed] (bU) -- (bW);
\draw[removed] (bV) -- (bW);
\draw[added] (aU) -- (aP);
\draw[added] (aW) -- (aP);
\draw[added] (aV) -- (aQ);
\draw[added] (aW) -- (aQ);
\draw[added] (aP) -- (aQ);
\end{scope}
\draw[->,thick] (1.4,0)--(2.6,0);
\end{scope}
\begin{scope}[shift={(0,-5.4)}]
\node[anchor=east] at (-1.6,0) {$r=8$};
\node[oldpair] (bU) at (0,0.8) {$U$};
\node[oldpair] (bV) at (-0.8,-0.55) {$V$};
\node[oldpair] (bW) at (0.8,-0.55) {$W$};
\node[oldpair] (aU) at (4.3,0.44) {$U$};
\node[oldpair] (aV) at (4.3,0.0) {$V$};
\node[oldpair] (aW) at (4.3,-0.44) {$W$};
\node[newpair] (aP) at (6.8,.8) {$P$};
\node[newpair] (aQ) at (6.8,-.8) {$Q$};
\begin{scope}[on background layer]
\draw[removed] (bU) -- (bV);
\draw[removed] (bV) -- (bW);
\draw[removed] (bW) -- (bU);
\draw[added] (aU) -- (aP);
\draw[added] (aU) -- (aQ);
\draw[added] (aV) -- (aP);
\draw[added] (aV) -- (aQ);
\draw[added] (aW) -- (aP);
\draw[added] (aW) -- (aQ);
\draw[added] (aP) -- (aQ);
\end{scope}
\draw[->,thick] (1.4,0)--(2.6,0);
\end{scope}
\begin{scope}[shift={(0,-8.1)}]
\node[anchor=east] at (-1.6,0) {$r=10$};
\node[oldpair] (bU) at (-0.7,0.65) {$U$};
\node[oldpair] (bV) at (0.7,0.65) {$V$};
\node[oldpair] (bW) at (0.7,-0.65) {$W$};
\node[oldpair] (bX) at (-0.7,-0.65) {$X$};
\node[oldpair] (aU) at (4.3,0.66) {$U$};
\node[oldpair] (aV) at (4.3,0.22) {$V$};
\node[oldpair] (aW) at (4.3,-0.22) {$W$};
\node[oldpair] (aX) at (4.3,-0.66) {$X$};
\node[newpair] (aP) at (6.8,.8) {$P$};
\node[newpair] (aQ) at (6.8,-.8) {$Q$};
\begin{scope}[on background layer]
\draw[removed] (bU) -- (bV);
\draw[removed] (bV) -- (bW);
\draw[removed] (bW) -- (bX);
\draw[removed] (bX) -- (bU);
\draw[added] (aU) -- (aP);
\draw[added] (aU) -- (aQ);
\draw[added] (aV) -- (aP);
\draw[added] (aV) -- (aQ);
\draw[added] (aW) -- (aP);
\draw[added] (aW) -- (aQ);
\draw[added] (aX) -- (aP);
\draw[added] (aX) -- (aQ);
\draw[added] (aP) -- (aQ);
\end{scope}
\draw[->,thick] (1.4,0)--(2.6,0);
\end{scope}
\begin{scope}[shift={(0,-10.8)}]
\node[anchor=east] at (-1.6,0) {$r=12$};
\node[oldpair] (bU) at (0,0.85) {$U$};
\node[oldpair] (bV) at (0.85,0.2) {$V$};
\node[oldpair] (bX) at (0.55,-0.7) {$X$};
\node[oldpair] (bW) at (-0.55,-0.7) {$W$};
\node[oldpair] (bY) at (-0.85,0.2) {$Y$};
\node[oldpair] (aU) at (4.3,0.88) {$U$};
\node[oldpair] (aV) at (4.3,0.44) {$V$};
\node[oldpair] (aX) at (4.3,0.0) {$X$};
\node[oldpair] (aW) at (4.3,-0.44) {$W$};
\node[oldpair] (aY) at (4.3,-0.88) {$Y$};
\node[newpair] (aP) at (6.8,.8) {$P$};
\node[newpair] (aQ) at (6.8,-.8) {$Q$};
\begin{scope}[on background layer]
\draw[removed] (bU) -- (bV);
\draw[removed] (bV) -- (bX);
\draw[removed] (bX) -- (bW);
\draw[removed] (bW) -- (bY);
\draw[removed] (bY) -- (bU);
\draw[retained] (bU) -- (bX);
\draw[retained] (aU) to[out=180,in=180,looseness=1.6] (aX);
\draw[added] (aU) -- (aP);
\draw[added] (aU) -- (aQ);
\draw[added] (aV) -- (aP);
\draw[added] (aV) -- (aQ);
\draw[added] (aW) -- (aP);
\draw[added] (aW) -- (aQ);
\draw[added] (aX) -- (aP);
\draw[added] (aX) -- (aQ);
\draw[added] (aY) -- (aP);
\draw[added] (aY) -- (aQ);
\draw[added] (aP) -- (aQ);
\end{scope}
\draw[->,thick] (1.4,0)--(2.6,0);
\end{scope}
\begin{scope}[shift={(0,-13.5)}]
\node[anchor=east] at (-1.6,0) {$r=14$};
\node[oldpair] (bU) at (0,0.95) {$U$};
\node[oldpair] (bX) at (-0.75,0.3) {$X$};
\node[oldpair] (bW) at (0.75,0.3) {$W$};
\node[oldpair] (bY) at (-0.75,-0.3) {$Y$};
\node[oldpair] (bZ) at (0.75,-0.3) {$Z$};
\node[oldpair] (bV) at (0,-0.95) {$V$};
\node[oldpair] (aU) at (4.3,1.1) {$U$};
\node[oldpair] (aX) at (4.3,0.66) {$X$};
\node[oldpair] (aW) at (4.3,0.22) {$W$};
\node[oldpair] (aY) at (4.3,-0.22) {$Y$};
\node[oldpair] (aZ) at (4.3,-0.66) {$Z$};
\node[oldpair] (aV) at (4.3,-1.1) {$V$};
\node[newpair] (aP) at (6.8,.8) {$P$};
\node[newpair] (aQ) at (6.8,-.8) {$Q$};
\begin{scope}[on background layer]
\draw[removed] (bU) -- (bX);
\draw[removed] (bU) -- (bW);
\draw[removed] (bX) -- (bW);
\draw[removed] (bV) -- (bY);
\draw[removed] (bV) -- (bZ);
\draw[removed] (bY) -- (bZ);
\draw[retained] (bX) -- (bY);
\draw[retained] (bW) -- (bZ);
\draw[retained] (aX) to[out=180,in=180,looseness=1.6] (aY);
\draw[retained] (aW) to[out=180,in=180,looseness=1.6] (aZ);
\draw[added] (aU) -- (aP);
\draw[added] (aU) -- (aQ);
\draw[added] (aV) -- (aP);
\draw[added] (aV) -- (aQ);
\draw[added] (aW) -- (aP);
\draw[added] (aW) -- (aQ);
\draw[added] (aX) -- (aP);
\draw[added] (aX) -- (aQ);
\draw[added] (aY) -- (aP);
\draw[added] (aY) -- (aQ);
\draw[added] (aZ) -- (aP);
\draw[added] (aZ) -- (aQ);
\draw[added] (aP) -- (aQ);
\end{scope}
\draw[->,thick] (1.4,0)--(2.6,0);
\end{scope}
\draw[removed] (-1,-15.1)--(-.3,-15.1);
\node[right,font=\small] at (-.3,-15.1) {deleted};
\draw[retained] (1.8,-15.1)--(2.5,-15.1);
\node[right,font=\small] at (2.5,-15.1) {retained};
\draw[added] (4.9,-15.1)--(5.6,-15.1);
\node[right,font=\small] at (5.6,-15.1) {added};
\end{tikzpicture}
\caption{The recursive edge exchanges of Lemma~\ref{lem:step}.
Each node is a prescribed vertex pair and each line is a complete
join $K_{2,2}$. The new pairs $P,Q$ receive disjoint unused inverse
pairs of labels (only $P$ is added for $r=4$).
Dashed red joins are deleted, solid blue joins are retained, and
solid green joins are added. Only joins selected for the operation
are shown. Every unshown edge of the ambient graph is unchanged.
The new special copy is selected from the right-hand diagram as
specified in the proof.}
\label{fig:recursive-exchanges}
\end{figure}
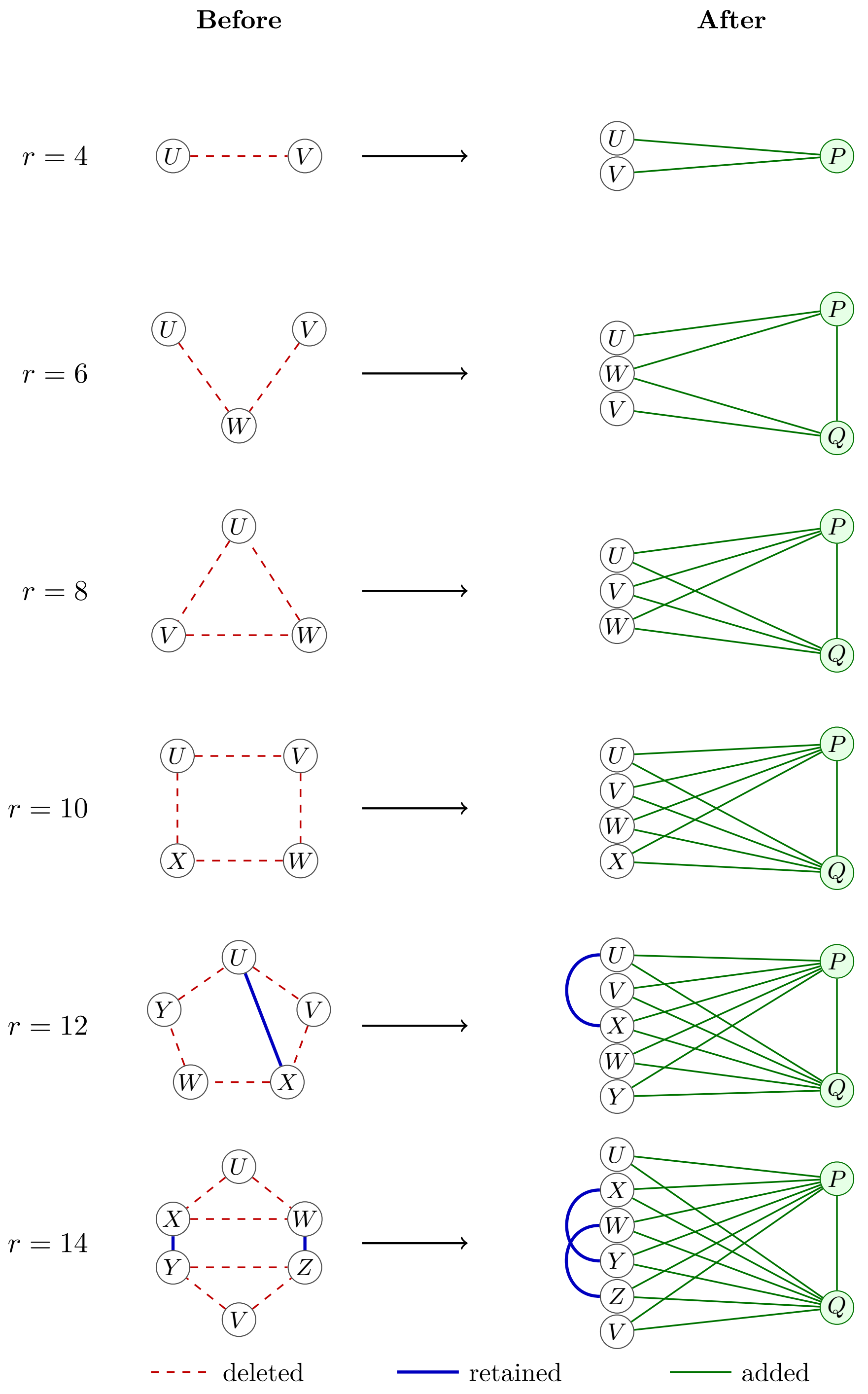

\begin{lemma}\label{lem:step}
Let $r\in\{4,6,8,10,12,14\}$, and let $G$ have an injective labeling
in $\Lambda$ with all weights zero and a special copy of $H_r$.
Suppose that the unused labels contain one two-element inverse pair
if $r=4$, and two disjoint such pairs if $r\ge6$.
There is an edge modification that adds vertices carrying exactly
these labels, with an injective labeling extending the original one.
It preserves all old degrees and weights and gives each new vertex
degree $r$ and weight zero. The resulting graph contains
a special copy of $H_r$ on old vertices of the selected copy and
new vertices. Vertices outside the old copy keep their neighbourhoods.
The operation preserves connectedness.
For $r=14$, the two new pairs occupy the positions $X,Z$ in the
new copy, with the names of Table~\ref{tab:patterns}.
\end{lemma}
\begin{proof}
We use the edge exchanges from Lemma~2.1 of~\cite{KFK4reg}
for $r=4$, Lemmas~1--4 of~\cite{KSKK6810} for $r=6,8,10,12$,
respectively, and Lemma~2.2 of~\cite{KovKrb14} for $r=14$.
Here we adapt these constructions to injective labelings in $\Lambda$
with all vertex weights zero.
The edge exchanges are shown in Figure~\ref{fig:recursive-exchanges}.
Label the new vertex pairs by the prescribed unused inverse pairs.
For $r=4$, delete $UV$ and add a new pair $P$ joined to $U,V$.
Old degrees are unchanged, new degrees are $4$, and $UP$ is a new
special $H_4$.

For $r=6$, delete $UW,VW$ and add new pairs $P,Q$ with joins
\[
 UP,WP,VQ,WQ,PQ.
\]
The pairs $U,V$ each lose and gain one join, while $W$ loses and
gains two. Each new vertex has degree $6$. The joins $UP,WP$
give a new special $H_6$.

For $r\in\{8,10,12,14\}$, delete $R_r$, retain $S_r$, and add new
pairs $P,Q$, joined to each other and to every pair of $H_r$.
Every vertex of the selected copy of $H_r$ loses and gains four edges. Since $H_r$ has
$(r-2)/2$ pairs, every new vertex has degree
$2((r-2)/2+1)=r$.
The following selected joins give a new special copy of $H_r$.
\[
\begin{array}{c|l}
r&\text{new copy of }H_r\\\hline
8&\text{triangle }U,P,Q\\
10&\text{cycle }V,P,W,Q,V\\
12&\text{cycle }W,P,X,U,Q,W\text{ with chord }PQ\\
14&\text{triangles }P,X,Y\text{ and }Q,W,Z\text{ with joins }PW,YQ.
\end{array}
\]
For $r=14$, the selected copy agrees with Table~\ref{tab:patterns}
under the renaming
\[
 (U',V',X',W',Y',Z')=(X,Z,P,Y,W,Q).
\]
Here the primed letters denote the pair names in the new copy.
Thus $P,Q$ occupy the positions $X',Z'$. Each vertex of these pairs
has degree $6$ in the selected copy, and no join between them is
selected, although $PQ$ is present in the ambient graph.

All added edges have a new vertex as an endpoint, so the resulting graph is simple.
Lemma~\ref{lem:replacement} gives the weight assertions.
The new special copy uses only vertices of the old copy and new
vertices. Hence its vertices all have degree $r$ whenever those of
the old copy do.

In every case, the added joins form a connected subgraph containing
all vertices of the old copy and all new vertices. Every deleted
edge therefore has its endpoints joined by a path in the resulting
graph. Replacing deleted edges by these paths shows that any two
old vertices remain connected whenever they were connected before
the operation. Every new vertex has an old neighbour, so
connectedness is preserved.
\end{proof}

\section{Modified bases and final steps}
For $r\ge6$, each recursive step uses two inverse pairs.
When $|A|\equiv3\pmod4$, the number $n(|A|-1)/2$ of inverse pairs
outside $K$ is odd, so these steps alone cannot exhaust the label
group. We reserve one or three inverse pairs for a final operation.
For $r=6,10,14$, that operation also increases the degrees of two
specified old vertices. We therefore first construct a graph in
which precisely these vertices have degree less than $r$.
We call their vertex pair the \emph{deficient pair}.
Lemma~\ref{lem:modified-group-base} supplies these graphs with
$\Z_n$-distance magic labelings of magic constant $0$.
For the orders covered by Lemma~\ref{lem:mixed}, we modify graphs
admitting integer distance magic labelings to obtain graphs each with
a deficient pair and a centered labeling, then apply
Lemma~\ref{lem:reduction}. The remaining cases
$(r,n)=(6,13),(10,15),(14,19)$ are supplied by the computational
examples in Figure~\ref{fig:small-modified-bases}, obtained using
the MILP formulation in the Appendix.

For the modified-base constructions with $r=10,14$, the prescribed
copies of $H_{10}$ and $H_{14}$ in Figure~\ref{fig:Hr-schematic}
must be supplemented by additional vertex pairs and joins.
We first establish the existence of these configurations, shown in
the left-hand column of Figure~\ref{fig:mixed-bases}.

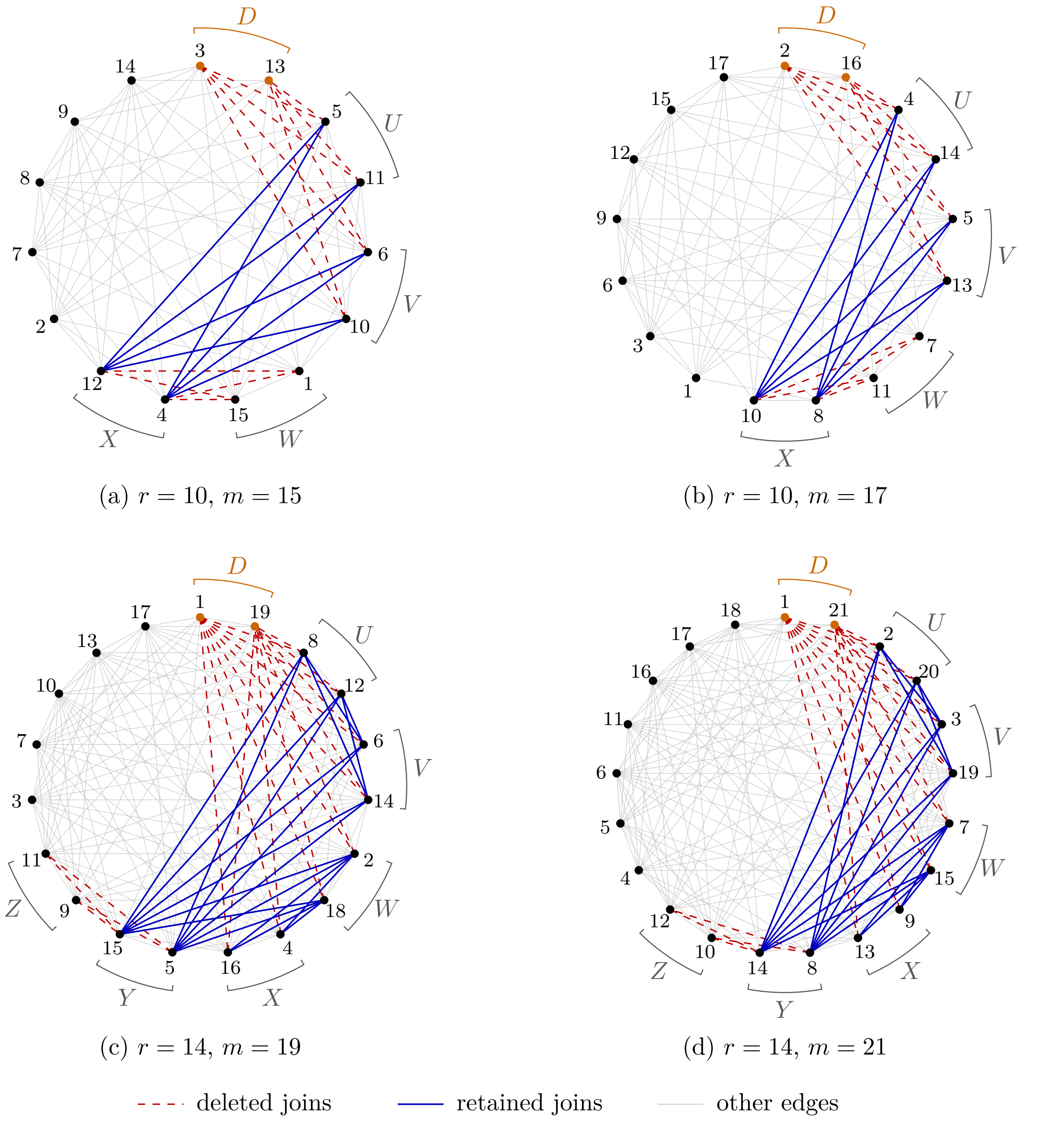
\begin{figure}[p]
\centering
\begin{minipage}{.48\textwidth}
\centering
\begin{tikzpicture}[
 x=1cm,y=1cm,
 seedvertex/.style={circle,draw=none,fill=black,
  minimum size=1.2mm,inner sep=0pt},
 seedother/.style={black!25,line width=.18pt},
 seedremoved/.style={red!75!black,dashed,line width=.55pt},
 seedretained/.style={blue!75!black,line width=.65pt}]
\path[use as bounding box] (-3.4,-3.4) rectangle (3.4,3.4);
\coordinate (v3) at (90.000000:2.4);
\coordinate (v13) at (66.000000:2.4);
\coordinate (v5) at (42.000000:2.4);
\coordinate (v11) at (18.000000:2.4);
\coordinate (v6) at (-6.000000:2.4);
\coordinate (v10) at (-30.000000:2.4);
\coordinate (v1) at (-54.000000:2.4);
\coordinate (v15) at (-78.000000:2.4);
\coordinate (v4) at (-102.000000:2.4);
\coordinate (v12) at (-126.000000:2.4);
\coordinate (v2) at (-150.000000:2.4);
\coordinate (v7) at (-174.000000:2.4);
\coordinate (v8) at (-198.000000:2.4);
\coordinate (v9) at (-222.000000:2.4);
\coordinate (v14) at (-246.000000:2.4);
\foreach \a/\b in {1/2,1/3,1/6,1/7,1/8,1/10,1/13,1/15,2/3,2/4,2/5,2/8,2/9,2/11,2/12,2/13,2/14,3/7,3/9,3/14,3/15,4/7,4/9,4/14,5/6,5/7,5/8,5/10,5/15,6/8,6/9,6/11,6/14,7/8,7/9,7/11,7/12,7/13,7/14,8/9,8/10,8/11,8/14,8/15,9/10,9/12,9/13,9/15,10/11,10/14,11/15,12/14,13/14,13/15,14/15}
  \draw[seedother] (v\a)--(v\b);
\foreach \a/\b in {1/4,1/12,3/5,3/6,3/10,3/11,4/15,5/13,6/13,10/13,11/13,12/15}
  \draw[seedremoved] (v\a)--(v\b);
\foreach \a/\b in {4/5,4/6,4/10,4/11,5/12,6/12,10/12,11/12}
  \draw[seedretained] (v\a)--(v\b);
\node[seedvertex,fill=orange!80!black] at (v3) {};
\node[font=\scriptsize,inner sep=0pt] at (90.000000:2.62) {$ 3 $};
\node[seedvertex,fill=orange!80!black] at (v13) {};
\node[font=\scriptsize,inner sep=0pt] at (66.000000:2.62) {$ 13 $};
\node[seedvertex] at (v5) {};
\node[font=\scriptsize,inner sep=0pt] at (42.000000:2.62) {$ 5 $};
\node[seedvertex] at (v11) {};
\node[font=\scriptsize,inner sep=0pt] at (18.000000:2.62) {$ 11 $};
\node[seedvertex] at (v6) {};
\node[font=\scriptsize,inner sep=0pt] at (-6.000000:2.62) {$ 6 $};
\node[seedvertex] at (v10) {};
\node[font=\scriptsize,inner sep=0pt] at (-30.000000:2.62) {$ 10 $};
\node[seedvertex] at (v1) {};
\node[font=\scriptsize,inner sep=0pt] at (-54.000000:2.62) {$ 1 $};
\node[seedvertex] at (v15) {};
\node[font=\scriptsize,inner sep=0pt] at (-78.000000:2.62) {$ 15 $};
\node[seedvertex] at (v4) {};
\node[font=\scriptsize,inner sep=0pt] at (-102.000000:2.62) {$ 4 $};
\node[seedvertex] at (v12) {};
\node[font=\scriptsize,inner sep=0pt] at (-126.000000:2.62) {$ 12 $};
\node[seedvertex] at (v2) {};
\node[font=\scriptsize,inner sep=0pt] at (-150.000000:2.62) {$ 2 $};
\node[seedvertex] at (v7) {};
\node[font=\scriptsize,inner sep=0pt] at (-174.000000:2.62) {$ 7 $};
\node[seedvertex] at (v8) {};
\node[font=\scriptsize,inner sep=0pt] at (-198.000000:2.62) {$ 8 $};
\node[seedvertex] at (v9) {};
\node[font=\scriptsize,inner sep=0pt] at (-222.000000:2.62) {$ 9 $};
\node[seedvertex] at (v14) {};
\node[font=\scriptsize,inner sep=0pt] at (-246.000000:2.62) {$ 14 $};
\draw[orange!80!black,line width=.45pt] (91.700000:2.86)--(91.700000:2.94) arc[start angle=91.700000,end angle=64.300000,radius=2.94cm]--(64.300000:2.86);
\node[font=\small,text=orange!80!black] at (78.000000:3.18) {$ D $};
\draw[black!65,line width=.45pt] (43.700000:2.86)--(43.700000:2.94) arc[start angle=43.700000,end angle=16.300000,radius=2.94cm]--(16.300000:2.86);
\node[font=\small,text=black!65] at (30.000000:3.18) {$ U $};
\draw[black!65,line width=.45pt] (-4.300000:2.86)--(-4.300000:2.94) arc[start angle=-4.300000,end angle=-31.700000,radius=2.94cm]--(-31.700000:2.86);
\node[font=\small,text=black!65] at (-18.000000:3.18) {$ V $};
\draw[black!65,line width=.45pt] (-52.300000:2.86)--(-52.300000:2.94) arc[start angle=-52.300000,end angle=-79.700000,radius=2.94cm]--(-79.700000:2.86);
\node[font=\small,text=black!65] at (-66.000000:3.18) {$ W $};
\draw[black!65,line width=.45pt] (-100.300000:2.86)--(-100.300000:2.94) arc[start angle=-100.300000,end angle=-127.700000,radius=2.94cm]--(-127.700000:2.86);
\node[font=\small,text=black!65] at (-114.000000:3.18) {$ X $};
\end{tikzpicture}
\par\small (a) $r=10$, $m=15$
\end{minipage}
\hfill
\begin{minipage}{.48\textwidth}
\centering
\begin{tikzpicture}[
 x=1cm,y=1cm,
 seedvertex/.style={circle,draw=none,fill=black,
  minimum size=1.2mm,inner sep=0pt},
 seedother/.style={black!25,line width=.18pt},
 seedremoved/.style={red!75!black,dashed,line width=.55pt},
 seedretained/.style={blue!75!black,line width=.65pt}]
\path[use as bounding box] (-3.4,-3.4) rectangle (3.4,3.4);
\coordinate (v2) at (90.000000:2.4);
\coordinate (v16) at (68.823529:2.4);
\coordinate (v4) at (47.647059:2.4);
\coordinate (v14) at (26.470588:2.4);
\coordinate (v5) at (5.294118:2.4);
\coordinate (v13) at (-15.882353:2.4);
\coordinate (v7) at (-37.058824:2.4);
\coordinate (v11) at (-58.235294:2.4);
\coordinate (v8) at (-79.411765:2.4);
\coordinate (v10) at (-100.588235:2.4);
\coordinate (v1) at (-121.764706:2.4);
\coordinate (v3) at (-142.941176:2.4);
\coordinate (v6) at (-164.117647:2.4);
\coordinate (v9) at (-185.294118:2.4);
\coordinate (v12) at (-206.470588:2.4);
\coordinate (v15) at (-227.647059:2.4);
\coordinate (v17) at (-248.823529:2.4);
\foreach \a/\b in {1/2,1/3,1/4,1/5,1/6,1/8,1/14,1/15,1/16,1/17,2/3,2/7,2/10,2/16,2/17,3/4,3/6,3/9,3/11,3/12,3/13,3/15,3/17,4/7,4/12,4/14,4/17,5/6,5/9,5/11,5/12,5/15,6/7,6/9,6/11,6/12,6/13,6/14,6/15,7/9,7/11,7/12,7/13,7/15,8/9,8/10,8/16,9/10,9/11,9/12,9/13,9/15,10/17,11/12,11/14,11/16,12/13,12/15,12/17,13/17,14/15,14/17,15/16,15/17,16/17}
  \draw[seedother] (v\a)--(v\b);
\foreach \a/\b in {2/4,2/5,2/13,2/14,4/16,5/16,7/8,7/10,8/11,10/11,13/16,14/16}
  \draw[seedremoved] (v\a)--(v\b);
\foreach \a/\b in {4/8,4/10,5/8,5/10,8/13,8/14,10/13,10/14}
  \draw[seedretained] (v\a)--(v\b);
\node[seedvertex,fill=orange!80!black] at (v2) {};
\node[font=\scriptsize,inner sep=0pt] at (90.000000:2.62) {$ 2 $};
\node[seedvertex,fill=orange!80!black] at (v16) {};
\node[font=\scriptsize,inner sep=0pt] at (68.823529:2.62) {$ 16 $};
\node[seedvertex] at (v4) {};
\node[font=\scriptsize,inner sep=0pt] at (47.647059:2.62) {$ 4 $};
\node[seedvertex] at (v14) {};
\node[font=\scriptsize,inner sep=0pt] at (26.470588:2.62) {$ 14 $};
\node[seedvertex] at (v5) {};
\node[font=\scriptsize,inner sep=0pt] at (5.294118:2.62) {$ 5 $};
\node[seedvertex] at (v13) {};
\node[font=\scriptsize,inner sep=0pt] at (-15.882353:2.62) {$ 13 $};
\node[seedvertex] at (v7) {};
\node[font=\scriptsize,inner sep=0pt] at (-37.058824:2.62) {$ 7 $};
\node[seedvertex] at (v11) {};
\node[font=\scriptsize,inner sep=0pt] at (-58.235294:2.62) {$ 11 $};
\node[seedvertex] at (v8) {};
\node[font=\scriptsize,inner sep=0pt] at (-79.411765:2.62) {$ 8 $};
\node[seedvertex] at (v10) {};
\node[font=\scriptsize,inner sep=0pt] at (-100.588235:2.62) {$ 10 $};
\node[seedvertex] at (v1) {};
\node[font=\scriptsize,inner sep=0pt] at (-121.764706:2.62) {$ 1 $};
\node[seedvertex] at (v3) {};
\node[font=\scriptsize,inner sep=0pt] at (-142.941176:2.62) {$ 3 $};
\node[seedvertex] at (v6) {};
\node[font=\scriptsize,inner sep=0pt] at (-164.117647:2.62) {$ 6 $};
\node[seedvertex] at (v9) {};
\node[font=\scriptsize,inner sep=0pt] at (-185.294118:2.62) {$ 9 $};
\node[seedvertex] at (v12) {};
\node[font=\scriptsize,inner sep=0pt] at (-206.470588:2.62) {$ 12 $};
\node[seedvertex] at (v15) {};
\node[font=\scriptsize,inner sep=0pt] at (-227.647059:2.62) {$ 15 $};
\node[seedvertex] at (v17) {};
\node[font=\scriptsize,inner sep=0pt] at (-248.823529:2.62) {$ 17 $};
\draw[orange!80!black,line width=.45pt] (91.700000:2.86)--(91.700000:2.94) arc[start angle=91.700000,end angle=67.123529,radius=2.94cm]--(67.123529:2.86);
\node[font=\small,text=orange!80!black] at (79.411765:3.18) {$ D $};
\draw[black!65,line width=.45pt] (49.347059:2.86)--(49.347059:2.94) arc[start angle=49.347059,end angle=24.770588,radius=2.94cm]--(24.770588:2.86);
\node[font=\small,text=black!65] at (37.058824:3.18) {$ U $};
\draw[black!65,line width=.45pt] (6.994118:2.86)--(6.994118:2.94) arc[start angle=6.994118,end angle=-17.582353,radius=2.94cm]--(-17.582353:2.86);
\node[font=\small,text=black!65] at (-5.294118:3.18) {$ V $};
\draw[black!65,line width=.45pt] (-35.358824:2.86)--(-35.358824:2.94) arc[start angle=-35.358824,end angle=-59.935294,radius=2.94cm]--(-59.935294:2.86);
\node[font=\small,text=black!65] at (-47.647059:3.18) {$ W $};
\draw[black!65,line width=.45pt] (-77.711765:2.86)--(-77.711765:2.94) arc[start angle=-77.711765,end angle=-102.288235,radius=2.94cm]--(-102.288235:2.86);
\node[font=\small,text=black!65] at (-90.000000:3.18) {$ X $};
\end{tikzpicture}
\par\small (b) $r=10$, $m=17$
\end{minipage}
\par\vspace{5mm}
\begin{minipage}{.48\textwidth}
\centering
\begin{tikzpicture}[
 x=1cm,y=1cm,
 seedvertex/.style={circle,draw=none,fill=black,
  minimum size=1.2mm,inner sep=0pt},
 seedother/.style={black!25,line width=.18pt},
 seedremoved/.style={red!75!black,dashed,line width=.55pt},
 seedretained/.style={blue!75!black,line width=.65pt}]
\path[use as bounding box] (-3.4,-3.4) rectangle (3.4,3.4);
\coordinate (v1) at (90.000000:2.4);
\coordinate (v19) at (71.052632:2.4);
\coordinate (v8) at (52.105263:2.4);
\coordinate (v12) at (33.157895:2.4);
\coordinate (v6) at (14.210526:2.4);
\coordinate (v14) at (-4.736842:2.4);
\coordinate (v2) at (-23.684211:2.4);
\coordinate (v18) at (-42.631579:2.4);
\coordinate (v4) at (-61.578947:2.4);
\coordinate (v16) at (-80.526316:2.4);
\coordinate (v5) at (-99.473684:2.4);
\coordinate (v15) at (-118.421053:2.4);
\coordinate (v9) at (-137.368421:2.4);
\coordinate (v11) at (-156.315789:2.4);
\coordinate (v3) at (-175.263158:2.4);
\coordinate (v7) at (-194.210526:2.4);
\coordinate (v10) at (-213.157895:2.4);
\coordinate (v13) at (-232.105263:2.4);
\coordinate (v17) at (-251.052632:2.4);
\foreach \a/\b in {1/3,1/5,1/7,1/9,1/17,1/19,2/3,2/7,2/8,2/9,2/10,2/11,2/14,2/18,3/4,3/6,3/8,3/9,3/10,3/11,3/12,3/13,3/14,3/15,3/17,3/18,4/5,4/6,4/7,4/10,4/11,4/12,4/13,4/16,4/17,5/10,5/13,5/15,5/17,6/7,6/9,6/11,6/13,6/16,6/17,7/8,7/9,7/10,7/11,7/12,7/13,7/14,7/15,7/16,7/19,8/10,8/12,8/13,8/16,8/17,9/10,9/11,9/13,9/14,9/16,9/18,9/19,10/11,10/12,10/13,10/15,10/16,10/17,10/18,11/13,11/14,11/16,11/17,11/19,12/13,12/17,12/18,13/14,13/16,13/17,13/18,14/17,14/18,15/16,15/19,16/17,17/18,17/19}
  \draw[seedother] (v\a)--(v\b);
\foreach \a/\b in {1/2,1/4,1/6,1/8,1/12,1/14,1/16,1/18,2/19,4/19,5/9,5/11,6/19,8/19,9/15,11/15,12/19,14/19,16/19,18/19}
  \draw[seedremoved] (v\a)--(v\b);
\foreach \a/\b in {2/4,2/5,2/15,2/16,4/18,5/6,5/8,5/12,5/14,5/18,6/8,6/12,6/15,8/14,8/15,12/14,12/15,14/15,15/18,16/18}
  \draw[seedretained] (v\a)--(v\b);
\node[seedvertex,fill=orange!80!black] at (v1) {};
\node[font=\scriptsize,inner sep=0pt] at (90.000000:2.62) {$ 1 $};
\node[seedvertex,fill=orange!80!black] at (v19) {};
\node[font=\scriptsize,inner sep=0pt] at (71.052632:2.62) {$ 19 $};
\node[seedvertex] at (v8) {};
\node[font=\scriptsize,inner sep=0pt] at (52.105263:2.62) {$ 8 $};
\node[seedvertex] at (v12) {};
\node[font=\scriptsize,inner sep=0pt] at (33.157895:2.62) {$ 12 $};
\node[seedvertex] at (v6) {};
\node[font=\scriptsize,inner sep=0pt] at (14.210526:2.62) {$ 6 $};
\node[seedvertex] at (v14) {};
\node[font=\scriptsize,inner sep=0pt] at (-4.736842:2.62) {$ 14 $};
\node[seedvertex] at (v2) {};
\node[font=\scriptsize,inner sep=0pt] at (-23.684211:2.62) {$ 2 $};
\node[seedvertex] at (v18) {};
\node[font=\scriptsize,inner sep=0pt] at (-42.631579:2.62) {$ 18 $};
\node[seedvertex] at (v4) {};
\node[font=\scriptsize,inner sep=0pt] at (-61.578947:2.62) {$ 4 $};
\node[seedvertex] at (v16) {};
\node[font=\scriptsize,inner sep=0pt] at (-80.526316:2.62) {$ 16 $};
\node[seedvertex] at (v5) {};
\node[font=\scriptsize,inner sep=0pt] at (-99.473684:2.62) {$ 5 $};
\node[seedvertex] at (v15) {};
\node[font=\scriptsize,inner sep=0pt] at (-118.421053:2.62) {$ 15 $};
\node[seedvertex] at (v9) {};
\node[font=\scriptsize,inner sep=0pt] at (-137.368421:2.62) {$ 9 $};
\node[seedvertex] at (v11) {};
\node[font=\scriptsize,inner sep=0pt] at (-156.315789:2.62) {$ 11 $};
\node[seedvertex] at (v3) {};
\node[font=\scriptsize,inner sep=0pt] at (-175.263158:2.62) {$ 3 $};
\node[seedvertex] at (v7) {};
\node[font=\scriptsize,inner sep=0pt] at (-194.210526:2.62) {$ 7 $};
\node[seedvertex] at (v10) {};
\node[font=\scriptsize,inner sep=0pt] at (-213.157895:2.62) {$ 10 $};
\node[seedvertex] at (v13) {};
\node[font=\scriptsize,inner sep=0pt] at (-232.105263:2.62) {$ 13 $};
\node[seedvertex] at (v17) {};
\node[font=\scriptsize,inner sep=0pt] at (-251.052632:2.62) {$ 17 $};
\draw[orange!80!black,line width=.45pt] (91.700000:2.86)--(91.700000:2.94) arc[start angle=91.700000,end angle=69.352632,radius=2.94cm]--(69.352632:2.86);
\node[font=\small,text=orange!80!black] at (80.526316:3.18) {$ D $};
\draw[black!65,line width=.45pt] (53.805263:2.86)--(53.805263:2.94) arc[start angle=53.805263,end angle=31.457895,radius=2.94cm]--(31.457895:2.86);
\node[font=\small,text=black!65] at (42.631579:3.18) {$ U $};
\draw[black!65,line width=.45pt] (15.910526:2.86)--(15.910526:2.94) arc[start angle=15.910526,end angle=-6.436842,radius=2.94cm]--(-6.436842:2.86);
\node[font=\small,text=black!65] at (4.736842:3.18) {$ V $};
\draw[black!65,line width=.45pt] (-21.984211:2.86)--(-21.984211:2.94) arc[start angle=-21.984211,end angle=-44.331579,radius=2.94cm]--(-44.331579:2.86);
\node[font=\small,text=black!65] at (-33.157895:3.18) {$ W $};
\draw[black!65,line width=.45pt] (-59.878947:2.86)--(-59.878947:2.94) arc[start angle=-59.878947,end angle=-82.226316,radius=2.94cm]--(-82.226316:2.86);
\node[font=\small,text=black!65] at (-71.052632:3.18) {$ X $};
\draw[black!65,line width=.45pt] (-97.773684:2.86)--(-97.773684:2.94) arc[start angle=-97.773684,end angle=-120.121053,radius=2.94cm]--(-120.121053:2.86);
\node[font=\small,text=black!65] at (-108.947368:3.18) {$ Y $};
\draw[black!65,line width=.45pt] (-135.668421:2.86)--(-135.668421:2.94) arc[start angle=-135.668421,end angle=-158.015789,radius=2.94cm]--(-158.015789:2.86);
\node[font=\small,text=black!65] at (-146.842105:3.18) {$ Z $};
\end{tikzpicture}
\par\small (c) $r=14$, $m=19$
\end{minipage}
\hfill
\begin{minipage}{.48\textwidth}
\centering
\begin{tikzpicture}[
 x=1cm,y=1cm,
 seedvertex/.style={circle,draw=none,fill=black,
  minimum size=1.2mm,inner sep=0pt},
 seedother/.style={black!25,line width=.18pt},
 seedremoved/.style={red!75!black,dashed,line width=.55pt},
 seedretained/.style={blue!75!black,line width=.65pt}]
\path[use as bounding box] (-3.4,-3.4) rectangle (3.4,3.4);
\coordinate (v1) at (90.000000:2.4);
\coordinate (v21) at (72.857143:2.4);
\coordinate (v2) at (55.714286:2.4);
\coordinate (v20) at (38.571429:2.4);
\coordinate (v3) at (21.428571:2.4);
\coordinate (v19) at (4.285714:2.4);
\coordinate (v7) at (-12.857143:2.4);
\coordinate (v15) at (-30.000000:2.4);
\coordinate (v9) at (-47.142857:2.4);
\coordinate (v13) at (-64.285714:2.4);
\coordinate (v8) at (-81.428571:2.4);
\coordinate (v14) at (-98.571429:2.4);
\coordinate (v10) at (-115.714286:2.4);
\coordinate (v12) at (-132.857143:2.4);
\coordinate (v4) at (-150.000000:2.4);
\coordinate (v5) at (-167.142857:2.4);
\coordinate (v6) at (-184.285714:2.4);
\coordinate (v11) at (-201.428571:2.4);
\coordinate (v16) at (-218.571429:2.4);
\coordinate (v17) at (-235.714286:2.4);
\coordinate (v18) at (-252.857143:2.4);
\foreach \a/\b in {1/4,1/5,1/10,1/12,1/17,1/18,2/5,2/6,2/7,2/10,2/11,2/15,2/16,2/18,3/4,3/6,3/9,3/11,3/12,3/13,3/16,3/17,4/5,4/6,4/7,4/8,4/10,4/11,4/14,4/15,4/16,4/18,4/19,4/21,5/6,5/8,5/9,5/11,5/12,5/13,5/14,5/16,5/17,5/20,5/21,6/7,6/9,6/10,6/12,6/13,6/15,6/17,6/18,6/19,6/20,7/11,7/12,7/16,7/17,7/20,8/9,8/13,8/17,8/18,9/10,9/11,9/14,9/16,9/18,9/19,10/11,10/12,10/13,10/16,10/17,10/20,10/21,11/12,11/13,11/15,11/17,11/18,11/19,11/20,12/15,12/16,12/18,12/19,12/21,13/14,13/16,13/18,13/19,14/17,14/18,15/16,15/17,15/20,16/17,16/18,16/19,16/20,17/18,17/19,17/21,18/20,18/21}
  \draw[seedother] (v\a)--(v\b);
\foreach \a/\b in {1/2,1/3,1/7,1/9,1/13,1/15,1/19,1/20,2/21,3/21,7/21,8/10,8/12,9/21,10/14,12/14,13/21,15/21,19/21,20/21}
  \draw[seedremoved] (v\a)--(v\b);
\foreach \a/\b in {2/3,2/8,2/14,2/19,3/8,3/14,3/20,7/8,7/9,7/13,7/14,8/15,8/19,8/20,9/15,13/15,14/15,14/19,14/20,19/20}
  \draw[seedretained] (v\a)--(v\b);
\node[seedvertex,fill=orange!80!black] at (v1) {};
\node[font=\scriptsize,inner sep=0pt] at (90.000000:2.62) {$ 1 $};
\node[seedvertex,fill=orange!80!black] at (v21) {};
\node[font=\scriptsize,inner sep=0pt] at (72.857143:2.62) {$ 21 $};
\node[seedvertex] at (v2) {};
\node[font=\scriptsize,inner sep=0pt] at (55.714286:2.62) {$ 2 $};
\node[seedvertex] at (v20) {};
\node[font=\scriptsize,inner sep=0pt] at (38.571429:2.62) {$ 20 $};
\node[seedvertex] at (v3) {};
\node[font=\scriptsize,inner sep=0pt] at (21.428571:2.62) {$ 3 $};
\node[seedvertex] at (v19) {};
\node[font=\scriptsize,inner sep=0pt] at (4.285714:2.62) {$ 19 $};
\node[seedvertex] at (v7) {};
\node[font=\scriptsize,inner sep=0pt] at (-12.857143:2.62) {$ 7 $};
\node[seedvertex] at (v15) {};
\node[font=\scriptsize,inner sep=0pt] at (-30.000000:2.62) {$ 15 $};
\node[seedvertex] at (v9) {};
\node[font=\scriptsize,inner sep=0pt] at (-47.142857:2.62) {$ 9 $};
\node[seedvertex] at (v13) {};
\node[font=\scriptsize,inner sep=0pt] at (-64.285714:2.62) {$ 13 $};
\node[seedvertex] at (v8) {};
\node[font=\scriptsize,inner sep=0pt] at (-81.428571:2.62) {$ 8 $};
\node[seedvertex] at (v14) {};
\node[font=\scriptsize,inner sep=0pt] at (-98.571429:2.62) {$ 14 $};
\node[seedvertex] at (v10) {};
\node[font=\scriptsize,inner sep=0pt] at (-115.714286:2.62) {$ 10 $};
\node[seedvertex] at (v12) {};
\node[font=\scriptsize,inner sep=0pt] at (-132.857143:2.62) {$ 12 $};
\node[seedvertex] at (v4) {};
\node[font=\scriptsize,inner sep=0pt] at (-150.000000:2.62) {$ 4 $};
\node[seedvertex] at (v5) {};
\node[font=\scriptsize,inner sep=0pt] at (-167.142857:2.62) {$ 5 $};
\node[seedvertex] at (v6) {};
\node[font=\scriptsize,inner sep=0pt] at (-184.285714:2.62) {$ 6 $};
\node[seedvertex] at (v11) {};
\node[font=\scriptsize,inner sep=0pt] at (-201.428571:2.62) {$ 11 $};
\node[seedvertex] at (v16) {};
\node[font=\scriptsize,inner sep=0pt] at (-218.571429:2.62) {$ 16 $};
\node[seedvertex] at (v17) {};
\node[font=\scriptsize,inner sep=0pt] at (-235.714286:2.62) {$ 17 $};
\node[seedvertex] at (v18) {};
\node[font=\scriptsize,inner sep=0pt] at (-252.857143:2.62) {$ 18 $};
\draw[orange!80!black,line width=.45pt] (91.700000:2.86)--(91.700000:2.94) arc[start angle=91.700000,end angle=71.157143,radius=2.94cm]--(71.157143:2.86);
\node[font=\small,text=orange!80!black] at (81.428571:3.18) {$ D $};
\draw[black!65,line width=.45pt] (57.414286:2.86)--(57.414286:2.94) arc[start angle=57.414286,end angle=36.871429,radius=2.94cm]--(36.871429:2.86);
\node[font=\small,text=black!65] at (47.142857:3.18) {$ U $};
\draw[black!65,line width=.45pt] (23.128571:2.86)--(23.128571:2.94) arc[start angle=23.128571,end angle=2.585714,radius=2.94cm]--(2.585714:2.86);
\node[font=\small,text=black!65] at (12.857143:3.18) {$ V $};
\draw[black!65,line width=.45pt] (-11.157143:2.86)--(-11.157143:2.94) arc[start angle=-11.157143,end angle=-31.700000,radius=2.94cm]--(-31.700000:2.86);
\node[font=\small,text=black!65] at (-21.428571:3.18) {$ W $};
\draw[black!65,line width=.45pt] (-45.442857:2.86)--(-45.442857:2.94) arc[start angle=-45.442857,end angle=-65.985714,radius=2.94cm]--(-65.985714:2.86);
\node[font=\small,text=black!65] at (-55.714286:3.18) {$ X $};
\draw[black!65,line width=.45pt] (-79.728571:2.86)--(-79.728571:2.94) arc[start angle=-79.728571,end angle=-100.271429,radius=2.94cm]--(-100.271429:2.86);
\node[font=\small,text=black!65] at (-90.000000:3.18) {$ Y $};
\draw[black!65,line width=.45pt] (-114.014286:2.86)--(-114.014286:2.94) arc[start angle=-114.014286,end angle=-134.557143,radius=2.94cm]--(-134.557143:2.86);
\node[font=\small,text=black!65] at (-124.285714:3.18) {$ Z $};
\end{tikzpicture}
\par\small (d) $r=14$, $m=21$
\end{minipage}
\par\vspace{3mm}
\begin{tikzpicture}[font=\small]
\draw[red!75!black,dashed,line width=.55pt] (0,0)--(.65,0);
\node[right] at (.7,0) {deleted joins};
\draw[blue!75!black,line width=.65pt] (3.7,0)--(4.35,0);
\node[right] at (4.4,0) {retained joins};
\draw[black!25,line width=.3pt] (7.4,0)--(8.05,0);
\node[right] at (8.1,0) {other edges};
\end{tikzpicture}
\caption{Initial regular graphs in Lemma~\ref{lem:initial-configurations},
used for the modified bases in Lemma~\ref{lem:mixed}.
Panels (a) and (b) redraw Figures~8 and~9
of~\cite{KSKK6810}, and panels (c) and (d) redraw Figures~4 and~5
of~\cite{KovKrb14}, with the original vertex labels.
Outer brackets identify the prescribed pairs, each of label sum $m+1$.
The pair $D$ is marked in orange. Dashed red edges form the joins
deleted in the proof, and solid blue edges form the retained joins
used to select the new copy of $H_r$. All remaining edges are gray.
Each panel shows the entire initial graph of order $m=n-2$.}
\label{fig:mixed-base-seeds}
\end{figure}

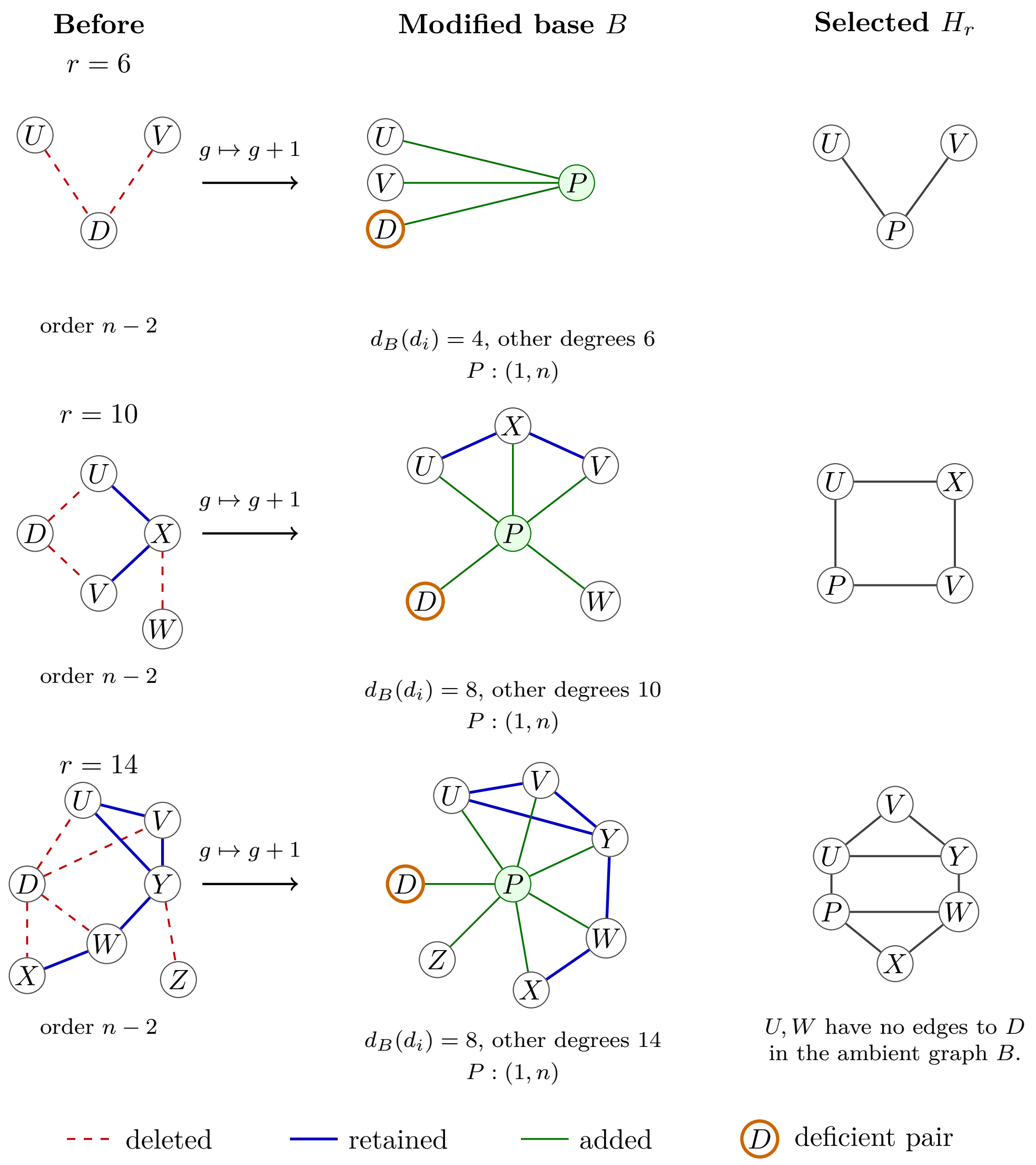
\begin{figure}[!htbp]
\centering
\begin{tikzpicture}[
 x=1cm,y=1cm,
 pair/.style={circle,draw=black!70,fill=white,minimum size=4.5mm,
  inner sep=.5pt,font=\small},
 newpair/.style={pair,draw=green!45!black,fill=green!10},
 deficient/.style={draw=orange!80!black,line width=1.2pt},
 removed/.style={red!75!black,dashed,line width=.7pt},
 retained/.style={blue!75!black,line width=1pt},
 added/.style={green!45!black,line width=.65pt},
 selected/.style={black!75,line width=.75pt}]
\node[font=\small\bfseries] at (0,2) {Before};
\node[font=\small\bfseries] at (5.2,2) {Modified base $B$};
\node[font=\small\bfseries] at (10,2) {Selected $H_r$};

\begin{scope}[shift={(0,-0.0)}]
\node[font=\small\bfseries] at (0,1.5) {$r=6$};
\node[pair] (bU) at (-0.8,0.6) {$U$};
\node[pair] (bV) at (0.8,0.6) {$V$};
\node[pair] (bD) at (0,-0.6) {$D$};
\node[pair] (aU) at (3.6,0.58) {$U$};
\node[pair] (aV) at (3.6,0.0) {$V$};
\node[pair,deficient] (aD) at (3.6,-0.58) {$D$};
\node[newpair] (aP) at (6,0) {$P$};
\node[pair] (sU) at (9.2,0.5) {$U$};
\node[pair] (sV) at (10.8,0.5) {$V$};
\node[pair] (sP) at (10,-0.6) {$P$};
\begin{scope}[on background layer]
\draw[removed] (bU) -- (bD);
\draw[removed] (bV) -- (bD);
\draw[added] (aU) -- (aP);
\draw[added] (aV) -- (aP);
\draw[added] (aD) -- (aP);
\draw[selected] (sU) -- (sP);
\draw[selected] (sV) -- (sP);
\end{scope}
\draw[->,thick] (1.3,0)--(2.5,0);
\node[above,font=\scriptsize] at (1.9,.15) {$g\mapsto g+1$};
\node[below,font=\scriptsize] at (0,-1.55) {order $n-2$};
\node[below,font=\scriptsize] at (5.2,-1.7) {$d_B(d_i)=4$, other degrees $6$};
\node[below,font=\scriptsize] at (5.2,-2.1) {$P:(1,n)$};
\end{scope}
\begin{scope}[shift={(0,-4.4)}]
\node[font=\small\bfseries] at (0,1.5) {$r=10$};
\node[pair] (bD) at (-0.8,0) {$D$};
\node[pair] (bU) at (0,0.75) {$U$};
\node[pair] (bV) at (0,-0.75) {$V$};
\node[pair] (bX) at (0.8,0) {$X$};
\node[pair] (bW) at (0.8,-1.2) {$W$};
\node[pair] (aU) at (4.1,0.85) {$U$};
\node[pair] (aV) at (6.3,0.85) {$V$};
\node[pair] (aX) at (5.2,1.35) {$X$};
\node[pair] (aW) at (6.3,-0.85) {$W$};
\node[pair,deficient] (aD) at (4.1,-0.85) {$D$};
\node[newpair] (aP) at (5.2,0) {$P$};
\node[pair] (sU) at (9.25,0.65) {$U$};
\node[pair] (sX) at (10.75,0.65) {$X$};
\node[pair] (sV) at (10.75,-0.65) {$V$};
\node[pair] (sP) at (9.25,-0.65) {$P$};
\begin{scope}[on background layer]
\draw[removed] (bD) -- (bU);
\draw[removed] (bD) -- (bV);
\draw[removed] (bW) -- (bX);
\draw[retained] (bU) -- (bX);
\draw[retained] (bV) -- (bX);
\draw[retained] (aU) -- (aX);
\draw[retained] (aV) -- (aX);
\draw[added] (aP) -- (aD);
\draw[added] (aP) -- (aU);
\draw[added] (aP) -- (aV);
\draw[added] (aP) -- (aW);
\draw[added] (aP) -- (aX);
\draw[selected] (sU) -- (sX);
\draw[selected] (sX) -- (sV);
\draw[selected] (sV) -- (sP);
\draw[selected] (sP) -- (sU);
\end{scope}
\draw[->,thick] (1.3,0)--(2.5,0);
\node[above,font=\scriptsize] at (1.9,.15) {$g\mapsto g+1$};
\node[below,font=\scriptsize] at (0,-1.55) {order $n-2$};
\node[below,font=\scriptsize] at (5.2,-1.7) {$d_B(d_i)=8$, other degrees $10$};
\node[below,font=\scriptsize] at (5.2,-2.1) {$P:(1,n)$};
\end{scope}
\begin{scope}[shift={(0,-8.8)}]
\node[font=\small\bfseries] at (0,1.5) {$r=14$};
\node[pair] (bD) at (-0.9,0) {$D$};
\node[pair] (bU) at (-0.2,1.05) {$U$};
\node[pair] (bV) at (0.8,0.8) {$V$};
\node[pair] (bY) at (0.8,0) {$Y$};
\node[pair] (bW) at (0.1,-0.75) {$W$};
\node[pair] (bX) at (-0.9,-1.15) {$X$};
\node[pair] (bZ) at (1,-1.2) {$Z$};
\node[pair,deficient] (aD) at (3.85,0) {$D$};
\node[pair] (aU) at (4.43,1.11) {$U$};
\node[pair] (aV) at (5.55,1.30) {$V$};
\node[pair] (aY) at (6.42,0.57) {$Y$};
\node[pair] (aW) at (6.37,-0.68) {$W$};
\node[pair] (aX) at (5.43,-1.33) {$X$};
\node[pair] (aZ) at (4.25,-0.95) {$Z$};
\node[newpair] (aP) at (5.2,0) {$P$};
\node[pair] (sV) at (10,1) {$V$};
\node[pair] (sU) at (9.2,0.35) {$U$};
\node[pair] (sY) at (10.8,0.35) {$Y$};
\node[pair] (sP) at (9.2,-0.35) {$P$};
\node[pair] (sW) at (10.8,-0.35) {$W$};
\node[pair] (sX) at (10,-1) {$X$};
\begin{scope}[on background layer]
\draw[removed] (bD) -- (bU);
\draw[removed] (bD) -- (bV);
\draw[removed] (bD) -- (bW);
\draw[removed] (bD) -- (bX);
\draw[removed] (bY) -- (bZ);
\draw[retained] (bU) -- (bV);
\draw[retained] (bW) -- (bX);
\draw[retained] (bY) -- (bU);
\draw[retained] (bY) -- (bV);
\draw[retained] (bY) -- (bW);
\draw[retained] (aU) -- (aV);
\draw[retained] (aW) -- (aX);
\draw[retained] (aY) -- (aU);
\draw[retained] (aY) -- (aV);
\draw[retained] (aY) -- (aW);
\draw[added] (aP) -- (aD);
\draw[added] (aP) -- (aU);
\draw[added] (aP) -- (aV);
\draw[added] (aP) -- (aW);
\draw[added] (aP) -- (aX);
\draw[added] (aP) -- (aY);
\draw[added] (aP) -- (aZ);
\draw[selected] (sU) -- (sV);
\draw[selected] (sV) -- (sY);
\draw[selected] (sY) -- (sU);
\draw[selected] (sW) -- (sX);
\draw[selected] (sX) -- (sP);
\draw[selected] (sP) -- (sW);
\draw[selected] (sU) -- (sP);
\draw[selected] (sY) -- (sW);
\end{scope}
\draw[->,thick] (1.3,0)--(2.5,0);
\node[above,font=\scriptsize] at (1.9,.15) {$g\mapsto g+1$};
\node[below,font=\scriptsize] at (0,-1.55) {order $n-2$};
\node[below,font=\scriptsize] at (5.2,-1.7) {$d_B(d_i)=8$, other degrees $14$};
\node[below,font=\scriptsize] at (5.2,-2.1) {$P:(1,n)$};
\node[below,align=center,font=\scriptsize] at (10,-1.55) {$U,W$ have no edges to $D$\\in the ambient graph $B$.};
\end{scope}
\draw[removed] (-.4,-12)--(.2,-12);
\node[right,font=\small] at (.2,-12) {deleted};
\draw[retained] (2.4,-12)--(3,-12);
\node[right,font=\small] at (3,-12) {retained};
\draw[added] (5.3,-12)--(5.9,-12);
\node[right,font=\small] at (5.9,-12) {added};
\node[pair,deficient] at (8.3,-12) {$D$};
\node[right,font=\small] at (8.6,-12) {deficient pair};
\end{tikzpicture}
\caption{The three modified-base constructions in Lemma~\ref{lem:mixed}.
Each node is a prescribed vertex pair and each line is a complete join
$K_{2,2}$. Green-filled nodes are new pairs. The orange outline marks
$D$, whose vertices have degree less than $r$ in $B$.
In each case, $D$ is an old pair and $P$ is the only new pair.
Only joins involved in the construction are displayed. All other
edges are unchanged. The right-hand column specifies the selected
copy of $H_r$ among vertices of degree $r$. Extra ambient edges are
omitted. After the shift, every prescribed pair has integer label sum
$n+1$. The new labels are listed below each modified base.}
\label{fig:mixed-bases}
\end{figure}

\Needspace{16\baselineskip}
\begin{lemma}\label{lem:initial-configurations}
Let $r\in\{10,14\}$ and let $m$ be an odd integer, with $m\ge15$
if $r=10$ and $m\ge19$ if $r=14$.
There is a connected $r$-regular graph of order $m$ with a distance
magic labeling $g$ and prescribed pairwise disjoint vertex pairs,
each of label sum $m+1$, satisfying the following conditions.
\begin{enumerate}
\item If $r=10$, the prescribed pairs are $D,U,V,W,X$ and the joins
\[
 DU,\quad DV,\quad WX,\quad UX,\quad VX
\]
are present.
\item If $r=14$, the prescribed pairs are $D,U,V,W,X,Y,Z$ and the joins
\[
 DU,\quad DV,\quad DW,\quad DX,\quad YZ,
 \qquad UV,\quad WX,\quad YU,\quad YV,\quad YW
\]
are present.
\end{enumerate}
\end{lemma}
\begin{proof}
Use the connected regular graphs and integer distance magic labelings
from the constructions in Proposition~\ref{prop:bases}. We verify
the required configurations for these graphs.
The recursive edge exchanges for $r=10,14$ are also shown in
Figure~\ref{fig:recursive-exchanges}. In the integer-labeled constructions,
a step from order $m-4$ to order $m$ increases every old label by $2$
and assigns the label pairs $(1,m)$ and $(2,m-1)$ to the two new
vertex pairs. Hence all prescribed old pairs and both new pairs have
label sum $m+1$. When referring to these right-hand diagrams, we
rename the pairs according to the present lemma, independently of
the names in the figure.

For $r=10$ and $m=15,17$, take the graphs in Figures~8 and~9
of~\cite{KSKK6810}, redrawn in Figure~\ref{fig:mixed-base-seeds}(a)
and~(b), respectively. Choose the pairs in the
following table. Each pair is identified by its two labels under $g$.
\[
\begin{array}{c|ccccc}
m&D&U&V&W&X\\\hline
15&(3,13)&(5,11)&(6,10)&(1,15)&(4,12)\\
17&(2,16)&(4,14)&(5,13)&(7,11)&(8,10)
\end{array}
\]
For every odd $m\ge19$, use the graph obtained by the construction
of Lemma~3 of~\cite{KSKK6810} from a graph of order $m-4$.
The edge exchange in this construction is shown in the $r=10$ row
of Figure~\ref{fig:recursive-exchanges}.
In the right-hand diagram, name the two new pairs $D,X$
and any three of the four old pairs $U,V,W$. Each of $D,X$ is
joined to all four old pairs, so the five required joins are present.
All these pairs have label sum $m+1$ under the resulting labeling $g$.

For $r=14$ and $m=19,21$, take the graphs in Figures~4 and~5
of~\cite{KovKrb14}, redrawn in Figure~\ref{fig:mixed-base-seeds}(c)
and~(d), respectively. Choose the pairs below, each identified
by its two labels under $g$.
\[
\begin{array}{c|ccccccc}
m&D&U&V&W&X&Y&Z\\\hline
19&(1,19)&(8,12)&(6,14)&(2,18)&(4,16)&(5,15)&(9,11)\\
21&(1,21)&(2,20)&(3,19)&(7,15)&(9,13)&(8,14)&(10,12)
\end{array}
\]
For every odd $m\ge23$, apply Lemma~2.2 of~\cite{KovKrb14}
to the corresponding graph of order $m-4$.
The edge exchange in this construction is shown in the $r=14$ row
of Figure~\ref{fig:recursive-exchanges}.
In the right-hand diagram, name the two new pairs $D,Y$. Name the endpoints of the two retained blue joins
$U,V$ and $W,X$, and let $Z$ be either of the two remaining old
pairs. Each of $D,Y$ is joined to all six old pairs, so all ten
required joins are present. Each chosen pair has label sum $m+1$.
Thus the required configuration exists for every odd $m\ge19$.

\end{proof}

\begin{lemma}\label{lem:mixed}
Let $n$ be a positive odd integer. For each row of the following
table, there is a connected graph $B$ of order $n$ admitting a
centered labeling and containing a specified pair $D=\{d_1,d_2\}$
whose vertices have degree $d$, while every vertex outside $D$ has
degree $r$. The graph contains a specified copy of $H_r$ disjoint
from $D$. Every prescribed pair of this copy, and $D$, has label
sum $n+1$.
\[
\begin{array}{c|c|c}
r&d&n\\\hline
6&4&n=11\text{ or }n\ge15\\
10&8&n\ge17\\
14&8&n\ge21
\end{array}
\]
For $r=14$, the copy can be named as in Table~\ref{tab:patterns}
so that no edge of $B$ joins a vertex of $X\cup Z$ to a vertex of $D$.
\end{lemma}
\begin{proof}
The constructions are shown in Figure~\ref{fig:mixed-bases}, including
the selected new copies of $H_r$.
Put $m=n-2$. For $r=6$, take the graph and integer labeling $g$
of order $m$ from Proposition~\ref{prop:bases}. For $r=10,14$, take
the graph, labeling $g$, and prescribed pairs supplied by
Lemma~\ref{lem:initial-configurations}.
For $r=6$, name its specified pairs $U,V,D$,
with joins $UD,VD$. Increase every old label by $1$, add a new pair
$P$ labeled $1,n$, delete $UD,VD$, and add $UP,VP,DP$.
Only the vertices of $D$ lose degree, which decreases from $6$ to $4$. The new vertices have degree $6$, and the joins $UP,VP$
form the required copy of $H_6$.
Before the edge changes, every old weight is $3(n+1)$.
At $D$, two neighbouring pairs are replaced by one, so the weight
becomes $2(n+1)$. All other old weights stay unchanged.
The new vertices see three pairs of sum $n+1$. Thus the labeling is centered.

For $r=10$, increase every old label by $1$ and add a new pair $P$ labeled $1,n$.
Delete $DU,DV,WX$, retain $UX,VX$, and add
\[
 PD,\quad PU,\quad PV,\quad PW,\quad PX.
\]
The vertices of $D$ lose two joins and gain one, so their degree
becomes $8$. Every vertex in $U\cup V\cup W\cup X$ loses and gains
one join, and each new vertex has degree $10$. All remaining
degrees are unchanged. The retained joins $UX,VX$ and the added
joins $PU,PV$ form a copy of $H_{10}$ on the cycle of pairs
$U,X,V,P,U$, disjoint from $D$.
The label shift gives every old vertex weight
$5(m+1)+10=5(n+1)$ and makes each of the five old pair sums equal
to $n+1$. The edge changes leave this weight unchanged outside $D$
and reduce it to $4(n+1)$ at $D$. Each new vertex has weight
$5(n+1)$, since its neighbourhood consists of five pairs of sum
$n+1$. The new labeling is a bijection onto $\{1,\ldots,n\}$ and
is centered.

For $r=14$, increase every old label by $1$ and add a new pair $P$ labeled $1,n$.
Delete $DU,DV,DW,DX,YZ$ and join $P$ to each of $D,U,V,W,X,Y,Z$.
Retain $UV,WX,YU,YV,YW$ and all other old edges.
Every vertex of $D$ loses four joins and gains one, so its degree
becomes $8$. Every vertex in $U\cup V\cup W\cup X\cup Y\cup Z$
loses and gains one join, and each new vertex has degree $14$.
All remaining degrees are unchanged.
The label shift gives every old vertex weight
$7(m+1)+14=7(n+1)$ and makes each chosen old pair sum equal to
$n+1$. The edge changes leave this weight unchanged outside $D$
and reduce it to $4(n+1)$ at $D$. Each new vertex has weight
$7(n+1)$, since its neighbourhood consists of seven pairs of sum
$n+1$. The new labeling is a bijection onto $\{1,\ldots,n\}$ and
is centered.

Select the triangles of pairs $U,V,Y$ and $W,X,P$, together with
joins $UP,YW$. They form a copy of $H_{14}$ disjoint from $D$.
Its pairs can be named as in Table~\ref{tab:patterns} by setting
\[
 (U',V',W',X',Y',Z')=(V,X,Y,U,P,W).
\]
The joins $DU,DW$ have been deleted, and no edge from $D$ to
$U\cup W$ has been added. Thus the pairs in the positions $X',Z'$
have no edges to $D$, as required.

The initial regular graph is connected by its choice.
In each case, the new pair $P$ is joined to every endpoint of a
deleted edge. Thus each deleted edge can be replaced by a path
through $P$, and every new vertex is connected to an old vertex.
The graph $B$ is therefore connected in all three cases.
\end{proof}

\Needspace{16\baselineskip}
\begin{lemma}\label{lem:modified-group-base}
Let $n$ be a positive odd integer. For each row of the following
table, there is a connected graph $B$ of order $n$ with a
$\Z_n$-distance magic labeling of magic constant $0$ and a
complementary pair $D$ whose vertices have degree $d$, while every
vertex outside $D$ has degree $r$.
\[
\begin{array}{c|c|c}
r&d&n\\\hline
6&4&n\ge11\\
10&8&n\ge15\\
14&8&n\ge19
\end{array}
\]
The graph contains a special copy of $H_r$ disjoint from $D$.
For $r=14$, its pairs can be named as in Table~\ref{tab:patterns}
so that no edge of $B$ joins a vertex of $X\cup Z$ to a vertex of $D$.
\end{lemma}
\begin{proof}
For the orders covered by Lemma~\ref{lem:mixed}, apply
Lemma~\ref{lem:reduction} to its centered labeling. This makes
$D$ and every prescribed pair complementary, while preserving
the graph and hence the required degrees and nonadjacencies.

The remaining cases are $(r,n)=(6,13),(10,15),(14,19)$.
Use the graphs in Figure~\ref{fig:small-modified-bases}(a), (b),
and (c), respectively, with the displayed residue labels and
prescribed pairs. Identifying each vertex with its label, direct
calculation from the displayed edge sets gives
\[
 d_B(v)=
 \begin{cases}
 d,&v\in D,\\
 r,&v\notin D,
 \end{cases}
 \qquad
 \sum_{u\in N_B(v)}u\equiv0\pmod n
\]
for every vertex $v$. The labels exhaust $\Z_n$, and each marked
pair has label sum $n$. The highlighted edges are precisely the
joins of $H_r$ in Table~\ref{tab:patterns}, and their endpoints
are disjoint from $D$. Thus these copies are special.
Every component contains a vertex outside $D$, since $d>1$ and
$|D|=2$, and therefore has at least $r+1$ vertices. In all three
cases $n<2(r+1)$, so the graphs are connected.
In panel (c), $X\cup Z=\{4,9,10,15\}$ and none of these vertices
is adjacent to either vertex of $D=\{1,18\}$. This verifies the
additional condition for $r=14$.
\end{proof}

\Needspace{20\baselineskip}
\begin{lemma}\label{lem:final}
Let $\Lambda$ be a finite Abelian group, let $r\in\{6,8,10,12,14\}$, and let $G$ have an injective labeling
$f\colon V(G)\to\Lambda$, with
all weights zero. Suppose that $G$ contains a special copy of
$H_r$ and that the following conditions hold.
\begin{enumerate}
\item If $r\in\{8,12\}$, then $G$ is $r$-regular.
\item If $r\in\{6,10,14\}$, then there is a complementary pair
$D=\{d_1,d_2\}$ disjoint from the copy of $H_r$, every vertex
outside $D$ has degree $r$, and
\[
d_G(d_1)=d_G(d_2)=
\begin{cases}
4,&r=6,\\
8,&r\in\{10,14\}.
\end{cases}
\]
\item If $r=14$, the pairs of the copy are named as in
Table~\ref{tab:patterns}, and no edge of $G$ joins a vertex of
$X\cup Z$ to a vertex of $D$.
\end{enumerate}
Set $s=1$ for $r\in\{6,10,14\}$ and $s=3$ for $r\in\{8,12\}$.
Suppose that the unused labels contain $s$ pairwise disjoint
two-element inverse pairs $L_1,\ldots,L_s$.
Then $G$ can be modified by adding $2s$ vertices and changing edges
to obtain an $r$-regular graph $G'$ with an injective labeling
$f'\colon V(G')\to\Lambda$ such that all weights are zero and
\[
 f'|_{V(G)}=f,\qquad
 f'\bigl(V(G')\setminus V(G)\bigr)=\bigcup_{i=1}^s L_i.
\]
If $G$ is connected, then $G'$ is connected.
\end{lemma}
\begin{figure}[p]
\centering
\begin{tikzpicture}[
 x=1.18cm,y=1.08cm,
 pair/.style={circle,draw=black!70,fill=white,minimum size=4.5mm,
  inner sep=.4pt,font=\footnotesize,text=black},
 newpair/.style={pair,draw=green!45!black,fill=green!10},
 distinguished/.style={draw=orange!80!black,line width=1pt},
 removed/.style={red!75!black,dashed,line width=.7pt},
 retained/.style={blue!75!black,line width=1.2pt},
 added/.style={green!45!black,line width=.6pt}]
\node[font=\small\bfseries] at (0,1.8) {Before};
\node[font=\small\bfseries] at (6,1.8) {After ($r$-regular)};

\begin{scope}[shift={(0,-0.0)}]
\node[anchor=east] at (-1.5,0) {$r=6$};
\node[pair] (bU) at (-0.8,0.6) {$U$};
\node[pair] (bV) at (0.8,0.6) {$V$};
\node[pair] (bW) at (0,-0.6) {$W$};
\node[pair,distinguished] (bD) at (1.4,-0.7) {$D$};
\node[pair] (aU) at (4.6,0.66) {$U$};
\node[pair] (aV) at (4.6,0.22) {$V$};
\node[pair] (aW) at (4.6,-0.22) {$W$};
\node[pair,distinguished] (aD) at (4.6,-0.66) {$D$};
\node[newpair] (aP) at (7,0) {$P$};
\begin{scope}[on background layer]
\draw[removed] (bU) -- (bW);
\draw[retained] (bV) -- (bW);
\draw[retained] (aV) to[out=180,in=180,looseness=1.4] (aW);
\draw[added] (aU) -- (aP);
\draw[added] (aW) -- (aP);
\draw[added] (aD) -- (aP);
\end{scope}
\draw[->,thick] (2.1,0)--(3.4,0);
\node[below,font=\scriptsize] at (2.7,-1.35) {degree of $d_i$\quad $4\longrightarrow 6$};
\end{scope}
\begin{scope}[shift={(0,-3.25)}]
\node[anchor=east] at (-1.5,0) {$r=8$};
\node[pair] (bU) at (0,0.8) {$U$};
\node[pair] (bV) at (-0.8,-0.6) {$V$};
\node[pair] (bW) at (0.8,-0.6) {$W$};
\node[pair] (aU) at (4.6,0.44) {$U$};
\node[pair] (aV) at (4.6,0) {$V$};
\node[pair] (aW) at (4.6,-0.44) {$W$};
\node[newpair] (aP) at (7,1) {$P$};
\node[newpair] (aQ) at (8,0) {$Q$};
\node[newpair] (aT) at (7,-1) {$T$};
\begin{scope}[on background layer]
\draw[removed] (bU) -- (bV);
\draw[removed] (bV) -- (bW);
\draw[removed] (bW) -- (bU);
\draw[added] (aP) -- (aQ);
\draw[added] (aP) -- (aT);
\draw[added] (aQ) -- (aT);
\draw[added] (aU) -- (aP);
\draw[added] (aU) -- (aQ);
\draw[added] (aV) -- (aP);
\draw[added] (aV) -- (aT);
\draw[added] (aW) -- (aQ);
\draw[added] (aW) -- (aT);
\end{scope}
\draw[->,thick] (2.1,0)--(3.4,0);
\end{scope}
\begin{scope}[shift={(0,-6.5)}]
\node[anchor=east] at (-1.5,0) {$r=10$};
\node[pair] (bU) at (-0.7,0.7) {$U$};
\node[pair] (bV) at (0.7,0.7) {$V$};
\node[pair] (bW) at (0.7,-0.7) {$W$};
\node[pair] (bX) at (-0.7,-0.7) {$X$};
\node[pair,distinguished] (bD) at (1.5,0) {$D$};
\node[pair] (aU) at (4.6,0.88) {$U$};
\node[pair] (aV) at (4.6,0.44) {$V$};
\node[pair] (aW) at (4.6,0.0) {$W$};
\node[pair] (aX) at (4.6,-0.44) {$X$};
\node[pair,distinguished] (aD) at (4.6,-0.88) {$D$};
\node[newpair] (aP) at (7,0) {$P$};
\begin{scope}[on background layer]
\draw[removed] (bU) -- (bV);
\draw[removed] (bW) -- (bX);
\draw[retained] (bV) -- (bW);
\draw[retained] (bX) -- (bU);
\draw[retained] (aV) to[out=180,in=180,looseness=1.4] (aW);
\draw[retained] (aX) to[out=180,in=180,looseness=1.4] (aU);
\draw[added] (aU) -- (aP);
\draw[added] (aV) -- (aP);
\draw[added] (aW) -- (aP);
\draw[added] (aX) -- (aP);
\draw[added] (aD) -- (aP);
\end{scope}
\draw[->,thick] (2.1,0)--(3.4,0);
\node[below,font=\scriptsize] at (2.7,-1.35) {degree of $d_i$\quad $8\longrightarrow 10$};
\end{scope}
\begin{scope}[shift={(0,-9.75)}]
\node[anchor=east] at (-1.5,0) {$r=12$};
\node[pair] (bU) at (0,0.9) {$U$};
\node[pair] (bV) at (0.85,0.2) {$V$};
\node[pair] (bX) at (0.55,-0.7) {$X$};
\node[pair] (bW) at (-0.55,-0.7) {$W$};
\node[pair] (bY) at (-0.85,0.2) {$Y$};
\node[pair] (aU) at (4.6,0.88) {$U$};
\node[pair] (aV) at (4.6,0.44) {$V$};
\node[pair] (aX) at (4.6,0.0) {$X$};
\node[pair] (aW) at (4.6,-0.44) {$W$};
\node[pair] (aY) at (4.6,-0.88) {$Y$};
\node[newpair] (aP) at (7,1) {$P$};
\node[newpair] (aQ) at (8,0) {$Q$};
\node[newpair] (aT) at (7,-1) {$T$};
\begin{scope}[on background layer]
\draw[removed] (bU) -- (bV);
\draw[removed] (bV) -- (bX);
\draw[removed] (bX) -- (bW);
\draw[removed] (bW) -- (bY);
\draw[removed] (bY) -- (bU);
\draw[removed] (bU) -- (bX);
\draw[added] (aP) -- (aQ);
\draw[added] (aP) -- (aT);
\draw[added] (aQ) -- (aT);
\draw[added] (aU) -- (aP);
\draw[added] (aU) -- (aQ);
\draw[added] (aU) -- (aT);
\draw[added] (aX) -- (aP);
\draw[added] (aX) -- (aQ);
\draw[added] (aX) -- (aT);
\draw[added] (aV) -- (aP);
\draw[added] (aV) -- (aQ);
\draw[added] (aW) -- (aQ);
\draw[added] (aW) -- (aT);
\draw[added] (aY) -- (aP);
\draw[added] (aY) -- (aT);
\end{scope}
\draw[->,thick] (2.1,0)--(3.4,0);
\end{scope}
\begin{scope}[shift={(0,-13.0)}]
\node[anchor=east] at (-1.5,0) {$r=14$};
\node[pair] (bU) at (0,1) {$U$};
\node[pair] (bX) at (-0.8,0.35) {$X$};
\node[pair] (bW) at (0.8,0.35) {$W$};
\node[pair] (bY) at (-0.8,-0.35) {$Y$};
\node[pair] (bZ) at (0.8,-0.35) {$Z$};
\node[pair] (bV) at (0,-1) {$V$};
\node[pair,distinguished] (bD) at (1.5,0) {$D$};
\node[pair] (aU) at (4.6,1.32) {$U$};
\node[pair] (aX) at (4.6,0.88) {$X$};
\node[pair] (aW) at (4.6,0.44) {$W$};
\node[pair] (aY) at (4.6,0.0) {$Y$};
\node[pair] (aZ) at (4.6,-0.44) {$Z$};
\node[pair] (aV) at (4.6,-0.88) {$V$};
\node[pair,distinguished] (aD) at (4.6,-1.32) {$D$};
\node[newpair] (aP) at (7,0) {$P$};
\begin{scope}[on background layer]
\draw[removed] (bU) -- (bX);
\draw[removed] (bX) -- (bW);
\draw[removed] (bV) -- (bZ);
\draw[removed] (bY) -- (bZ);
\draw[retained] (bU) -- (bW);
\draw[retained] (bV) -- (bY);
\draw[retained] (bX) -- (bY);
\draw[retained] (bW) -- (bZ);
\draw[retained] (aU) to[out=180,in=180,looseness=1.2] (aW);
\draw[retained] (aV) to[out=180,in=180,looseness=1.2] (aY);
\draw[retained] (aX) to[out=180,in=180,looseness=1.2] (aY);
\draw[retained] (aW) to[out=180,in=180,looseness=1.2] (aZ);
\draw[added] (aU) -- (aP);
\draw[added] (aV) -- (aP);
\draw[added] (aW) -- (aP);
\draw[added] (aX) -- (aP);
\draw[added] (aY) -- (aP);
\draw[added] (aZ) -- (aP);
\draw[added] (aD) -- (aP);
\draw[added] (aX) to[out=0,in=0,looseness=1.2] (aD);
\draw[added] (aZ) to[out=0,in=0,looseness=1.2] (aD);
\end{scope}
\draw[->,thick] (2.1,0)--(3.4,0);
\node[below,font=\scriptsize] at (2.7,-1.35) {degree of $d_i$\quad $8\longrightarrow 14$};
\end{scope}
\draw[removed] (-.4,-15.1)--(.2,-15.1);
\node[right,font=\small] at (.2,-15.1) {deleted};
\draw[retained] (2.6,-15.1)--(3.2,-15.1);
\node[right,font=\small] at (3.2,-15.1) {retained};
\draw[added] (5.7,-15.1)--(6.3,-15.1);
\node[right,font=\small] at (6.3,-15.1) {added};
\end{tikzpicture}
\caption{The final edge exchanges in Lemma~\ref{lem:final}.
Each node is a prescribed vertex pair and each line is a complete join
$K_{2,2}$. Green-filled pairs are new and receive disjoint unused
inverse pairs of labels. The orange outline identifies the old pair
$D$, whose degree increases to $r$.
Only joins used in the exchange are shown. All unshown edges remain
unchanged. In particular, $D$ may already have neighbours that are
not shown on the left. For $r=14$, the added old--old joins $DX,DZ$
are drawn as green arcs and must be absent before the exchange.
All retained joins are blue, including those outside $S_r$ when
$r=6,10,14$.}
\label{fig:final-exchanges}
\end{figure}
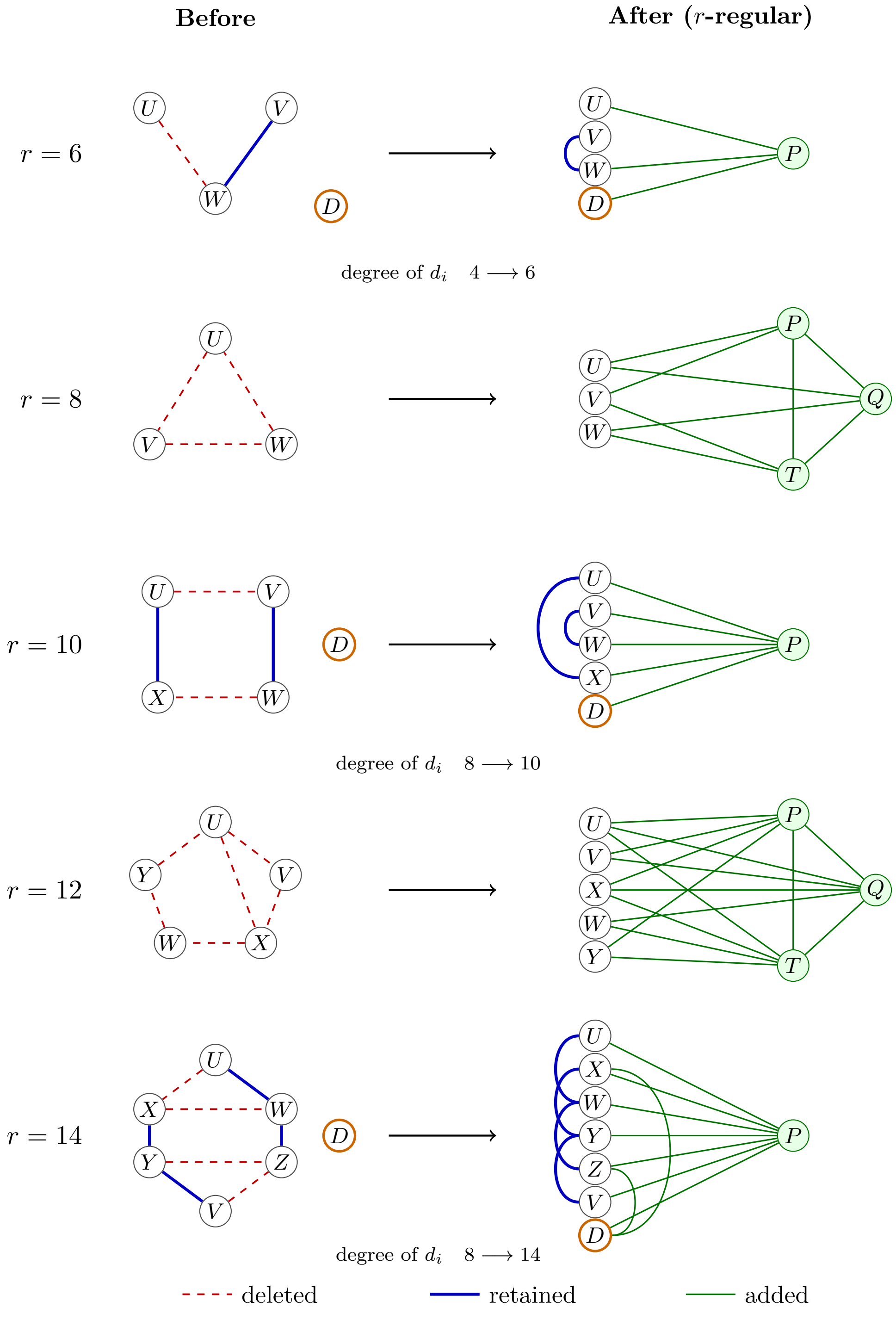

\begin{proof}
Figure~\ref{fig:final-exchanges} shows the five edge exchanges.
All new pairs receive pairwise disjoint unused two-element inverse pairs.
For $r=6$, delete $UW$ from the special $H_6$ and add a pair $P$
joined to $U,W,D$. The vertices of $U,W$ retain their degrees,
those of $D$ gain two neighbours, and each new vertex has degree $6$.
For $r=10$, delete $UV,WX$ and add a pair $P$ joined to $U,V,W,X,D$.
Again the specified old vertices retain their degrees, $D$ gains two
neighbours per vertex, and every new vertex has degree $10$.

For $r=14$, delete the joins $UX,XW,VZ,YZ$ and retain the other
joins of the specified $H_{14}$. Add one new pair $P$ joined to
each of $U,V,W,X,Y,Z,D$, and add the joins $DX,DZ$.
Each vertex in $U,V,W,Y$ loses two neighbours and gains the two
vertices of $P$. Each vertex in $X,Z$ loses four neighbours and
gains the four vertices of $P\cup D$. These vertices therefore
retain degree $14$. Each vertex of $D$ gains the six vertices of
$P\cup X\cup Z$, raising its degree from $8$ to $14$.
Each new vertex has the fourteen vertices of
$U\cup V\cup W\cup X\cup Y\cup Z\cup D$ as its neighbourhood.

For $r=8,12$, delete all joins of $H_r$, add three mutually joined
new pairs $P,Q,T$, and add the joins specified below.
\[
\begin{array}{c|l|l}
r&\text{old pair(s)}&\text{newly joined pairs}\\\hline
8&U& P,Q\\
 &V& P,T\\
 &W& Q,T\\\hline
12&U,X& P,Q,T\\
 &V& P,Q\\
 &W& Q,T\\
 &Y& P,T
\end{array}
\]
In both cases, each old vertex gains as many edges as it loses.
A new vertex has four neighbours in the other new pairs and,
respectively, four or eight old neighbours, giving degree $r$.

The only added edges between old vertices are $DX,DZ$, absent by
hypothesis. All other added edges have a new vertex as an endpoint, so the
resulting graphs are simple. Lemma~\ref{lem:replacement} proves
that all weights remain zero.

For $r=6,10,14$, the new pair is joined to both endpoints of every
deleted edge. For $r=8,12$, the mutually joined new pairs,
together with their joins to the old pairs, form a connected
subgraph containing all endpoints of deleted edges.
Every deleted edge can therefore be replaced by a path in $G'$,
and every new vertex has an old neighbour. Hence $G'$ is connected
whenever $G$ is connected.
\end{proof}

\FloatBarrier

\section{Proof of the main theorem}
Table~\ref{tab:construction-guide} summarizes the choices needed
when $|A|\equiv3\pmod4$ and $r\ge6$. In the cases with a deficient
pair $D$, the recursive steps keep it outside the selected copy of
$H_r$. Its vertices reach degree $r$ in the final step.

\begin{table}[htbp]
\centering
\caption{Construction guide for $|A|\equiv3\pmod4$ and $r\ge6$.
Every base has order $n$ and a $\Z_n$-distance magic labeling
with magic constant $0$.}
\label{tab:construction-guide}
\begin{tabular}{cccc}
\toprule
$r$ & Source of the base & \shortstack{Degree of each\\vertex in $D$} & Reserved pairs\\
\midrule
6 & Lemma~\ref{lem:modified-group-base} & 4 & 1\\
8 & Corollary~\ref{lem:special-base} & --- & 3\\
10 & Lemma~\ref{lem:modified-group-base} & 8 & 1\\
12 & Corollary~\ref{lem:special-base} & --- & 3\\
14 & Lemma~\ref{lem:modified-group-base} & 8 & 1\\
\bottomrule
\end{tabular}
\end{table}

Here ``---'' means that no deficient pair is needed. For $r=4$,
each recursive step uses a single inverse pair and no final step is
required.
\begin{figure}[!htbp]
\centering
\begin{minipage}[t]{.48\textwidth}
\centering
\begin{tikzpicture}[
  x=1cm,y=1cm,
  seedvertex/.style={circle,draw=none,fill=black,
    minimum size=1.2mm,inner sep=0pt},
  seedother/.style={black!30,line width=.25pt},
  seedremoved/.style={red!75!black,dashed,line width=.65pt}]
\path[use as bounding box] (-3.5,-3.45) rectangle (3.5,3.5);
\coordinate (v6) at (90.000000:2.4);
\coordinate (v7) at (62.307692:2.4);
\coordinate (v1) at (34.615385:2.4);
\coordinate (v12) at (6.923077:2.4);
\coordinate (v2) at (-20.769231:2.4);
\coordinate (v11) at (-48.461538:2.4);
\coordinate (v3) at (-76.153846:2.4);
\coordinate (v10) at (-103.846154:2.4);
\coordinate (v4) at (-131.538462:2.4);
\coordinate (v9) at (-159.230769:2.4);
\coordinate (v5) at (-186.923077:2.4);
\coordinate (v8) at (-214.615385:2.4);
\coordinate (v0) at (-242.307692:2.4);
\foreach \a/\b in {0/4,0/5,0/6,0/7,0/8,0/9,
  1/2,1/5,1/7,1/12,2/5,2/9,
  2/11,3/4,3/9,4/7,4/8,4/10,
  4/11,5/6,5/8,5/9,6/9,6/12,
  7/8,8/11,8/12,9/10,11/12}
  \draw[seedother] (v\a)--(v\b);
\foreach \a/\b in {1/3,1/10,2/3,2/10,3/11,3/12,10/11,10/12}
  \draw[seedremoved] (v\a)--(v\b);
\node[seedvertex,fill=orange!80!black] at (v6) {};
\node[font=\small,inner sep=0pt] at (90.000000:2.64) {$6$};
\node[seedvertex,fill=orange!80!black] at (v7) {};
\node[font=\small,inner sep=0pt] at (62.307692:2.64) {$7$};
\node[seedvertex] at (v1) {};
\node[font=\small,inner sep=0pt] at (34.615385:2.64) {$1$};
\node[seedvertex] at (v12) {};
\node[font=\small,inner sep=0pt] at (6.923077:2.64) {$12$};
\node[seedvertex] at (v2) {};
\node[font=\small,inner sep=0pt] at (-20.769231:2.64) {$2$};
\node[seedvertex] at (v11) {};
\node[font=\small,inner sep=0pt] at (-48.461538:2.64) {$11$};
\node[seedvertex] at (v3) {};
\node[font=\small,inner sep=0pt] at (-76.153846:2.64) {$3$};
\node[seedvertex] at (v10) {};
\node[font=\small,inner sep=0pt] at (-103.846154:2.64) {$10$};
\node[seedvertex] at (v4) {};
\node[font=\small,inner sep=0pt] at (-131.538462:2.64) {$4$};
\node[seedvertex] at (v9) {};
\node[font=\small,inner sep=0pt] at (-159.230769:2.64) {$9$};
\node[seedvertex] at (v5) {};
\node[font=\small,inner sep=0pt] at (-186.923077:2.64) {$5$};
\node[seedvertex] at (v8) {};
\node[font=\small,inner sep=0pt] at (-214.615385:2.64) {$8$};
\node[seedvertex] at (v0) {};
\node[font=\small,inner sep=0pt] at (-242.307692:2.64) {$0$};
\draw[orange!80!black,line width=.45pt] (91.700000:2.86)--(91.700000:2.94)
  arc[start angle=91.700000,end angle=60.607692,radius=2.94cm]--(60.607692:2.86);
\node[font=\small,text=orange!80!black] at (76.153846:3.18) {$D$};
\draw[black!65,line width=.45pt] (36.315385:2.86)--(36.315385:2.94)
  arc[start angle=36.315385,end angle=5.223077,radius=2.94cm]--(5.223077:2.86);
\node[font=\small,text=black!65] at (20.769231:3.18) {$U$};
\draw[black!65,line width=.45pt] (-19.069231:2.86)--(-19.069231:2.94)
  arc[start angle=-19.069231,end angle=-50.161538,radius=2.94cm]--(-50.161538:2.86);
\node[font=\small,text=black!65] at (-34.615385:3.18) {$V$};
\draw[black!65,line width=.45pt] (-74.453846:2.86)--(-74.453846:2.94)
  arc[start angle=-74.453846,end angle=-105.546154,radius=2.94cm]--(-105.546154:2.86);
\node[font=\small,text=black!65] at (-90.000000:3.18) {$W$};
\end{tikzpicture}
\par\small (a) $r=6$, $n=13$
\end{minipage}
\hfill
\begin{minipage}[t]{.48\textwidth}
\centering
\begin{tikzpicture}[
  x=1cm,y=1cm,
  seedvertex/.style={circle,draw=none,fill=black,
    minimum size=1.2mm,inner sep=0pt},
  seedother/.style={black!30,line width=.25pt},
  seedremoved/.style={red!75!black,dashed,line width=.65pt}]
\path[use as bounding box] (-3.5,-3.45) rectangle (3.5,3.5);
\coordinate (v3) at (90.000000:2.4);
\coordinate (v12) at (66.000000:2.4);
\coordinate (v1) at (42.000000:2.4);
\coordinate (v14) at (18.000000:2.4);
\coordinate (v6) at (-6.000000:2.4);
\coordinate (v9) at (-30.000000:2.4);
\coordinate (v7) at (-54.000000:2.4);
\coordinate (v8) at (-78.000000:2.4);
\coordinate (v5) at (-102.000000:2.4);
\coordinate (v10) at (-126.000000:2.4);
\coordinate (v2) at (-150.000000:2.4);
\coordinate (v13) at (-174.000000:2.4);
\coordinate (v4) at (-198.000000:2.4);
\coordinate (v11) at (-222.000000:2.4);
\coordinate (v0) at (-246.000000:2.4);
\foreach \a/\b in {0/1,0/2,0/3,0/4,0/7,0/8,
  0/11,0/12,0/13,0/14,1/2,1/3,
  1/4,1/7,1/14,2/3,2/4,2/5,
  2/6,2/8,2/9,2/11,2/13,3/4,
  3/8,3/9,3/10,3/11,4/5,4/6,
  4/8,4/10,4/12,4/13,5/6,5/10,
  5/11,5/12,6/9,6/12,6/13,7/8,
  7/11,7/12,7/13,8/14,9/10,9/11,
  9/13,10/11,10/13,11/12,11/13,11/14,
  12/13,12/14,13/14}
  \draw[seedother] (v\a)--(v\b);
\foreach \a/\b in {1/5,1/6,1/9,1/10,5/7,5/8,
  5/14,6/7,6/8,6/14,7/9,7/10,
  8/9,8/10,9/14,10/14}
  \draw[seedremoved] (v\a)--(v\b);
\node[seedvertex,fill=orange!80!black] at (v3) {};
\node[font=\small,inner sep=0pt] at (90.000000:2.64) {$3$};
\node[seedvertex,fill=orange!80!black] at (v12) {};
\node[font=\small,inner sep=0pt] at (66.000000:2.64) {$12$};
\node[seedvertex] at (v1) {};
\node[font=\small,inner sep=0pt] at (42.000000:2.64) {$1$};
\node[seedvertex] at (v14) {};
\node[font=\small,inner sep=0pt] at (18.000000:2.64) {$14$};
\node[seedvertex] at (v6) {};
\node[font=\small,inner sep=0pt] at (-6.000000:2.64) {$6$};
\node[seedvertex] at (v9) {};
\node[font=\small,inner sep=0pt] at (-30.000000:2.64) {$9$};
\node[seedvertex] at (v7) {};
\node[font=\small,inner sep=0pt] at (-54.000000:2.64) {$7$};
\node[seedvertex] at (v8) {};
\node[font=\small,inner sep=0pt] at (-78.000000:2.64) {$8$};
\node[seedvertex] at (v5) {};
\node[font=\small,inner sep=0pt] at (-102.000000:2.64) {$5$};
\node[seedvertex] at (v10) {};
\node[font=\small,inner sep=0pt] at (-126.000000:2.64) {$10$};
\node[seedvertex] at (v2) {};
\node[font=\small,inner sep=0pt] at (-150.000000:2.64) {$2$};
\node[seedvertex] at (v13) {};
\node[font=\small,inner sep=0pt] at (-174.000000:2.64) {$13$};
\node[seedvertex] at (v4) {};
\node[font=\small,inner sep=0pt] at (-198.000000:2.64) {$4$};
\node[seedvertex] at (v11) {};
\node[font=\small,inner sep=0pt] at (-222.000000:2.64) {$11$};
\node[seedvertex] at (v0) {};
\node[font=\small,inner sep=0pt] at (-246.000000:2.64) {$0$};
\draw[orange!80!black,line width=.45pt] (91.700000:2.86)--(91.700000:2.94)
  arc[start angle=91.700000,end angle=64.300000,radius=2.94cm]--(64.300000:2.86);
\node[font=\small,text=orange!80!black] at (78.000000:3.18) {$D$};
\draw[black!65,line width=.45pt] (43.700000:2.86)--(43.700000:2.94)
  arc[start angle=43.700000,end angle=16.300000,radius=2.94cm]--(16.300000:2.86);
\node[font=\small,text=black!65] at (30.000000:3.18) {$U$};
\draw[black!65,line width=.45pt] (-4.300000:2.86)--(-4.300000:2.94)
  arc[start angle=-4.300000,end angle=-31.700000,radius=2.94cm]--(-31.700000:2.86);
\node[font=\small,text=black!65] at (-18.000000:3.18) {$V$};
\draw[black!65,line width=.45pt] (-52.300000:2.86)--(-52.300000:2.94)
  arc[start angle=-52.300000,end angle=-79.700000,radius=2.94cm]--(-79.700000:2.86);
\node[font=\small,text=black!65] at (-66.000000:3.18) {$W$};
\draw[black!65,line width=.45pt] (-100.300000:2.86)--(-100.300000:2.94)
  arc[start angle=-100.300000,end angle=-127.700000,radius=2.94cm]--(-127.700000:2.86);
\node[font=\small,text=black!65] at (-114.000000:3.18) {$X$};
\end{tikzpicture}
\par\small (b) $r=10$, $n=15$
\end{minipage}
\par\medskip
\begin{minipage}[t]{.48\textwidth}
\centering
\begin{tikzpicture}[
  x=1cm,y=1cm,
  seedvertex/.style={circle,draw=none,fill=black,
    minimum size=1.2mm,inner sep=0pt},
  seedother/.style={black!30,line width=.25pt},
  seedremoved/.style={red!75!black,dashed,line width=.65pt},
  seedretained/.style={blue!75!black,line width=.8pt}]
\path[use as bounding box] (-3.5,-3.45) rectangle (3.5,3.5);
\coordinate (v1) at (90.000000:2.4);
\coordinate (v18) at (71.052632:2.4);
\coordinate (v7) at (52.105263:2.4);
\coordinate (v12) at (33.157895:2.4);
\coordinate (v9) at (14.210526:2.4);
\coordinate (v10) at (-4.736842:2.4);
\coordinate (v8) at (-23.684211:2.4);
\coordinate (v11) at (-42.631579:2.4);
\coordinate (v2) at (-61.578947:2.4);
\coordinate (v17) at (-80.526316:2.4);
\coordinate (v4) at (-99.473684:2.4);
\coordinate (v15) at (-118.421053:2.4);
\coordinate (v3) at (-137.368421:2.4);
\coordinate (v16) at (-156.315789:2.4);
\coordinate (v5) at (-175.263158:2.4);
\coordinate (v14) at (-194.210526:2.4);
\coordinate (v6) at (-213.157895:2.4);
\coordinate (v13) at (-232.105263:2.4);
\coordinate (v0) at (-251.052632:2.4);
\foreach \a/\b in {0/1,0/2,0/3,0/5,0/6,0/7,
  0/8,0/11,0/12,0/13,0/14,0/16,
  0/17,0/18,1/5,1/6,1/8,1/12,
  1/13,1/14,1/18,2/5,2/6,2/8,
  2/12,2/13,2/14,2/18,3/6,3/7,
  3/8,3/9,3/10,3/11,3/12,3/14,
  3/18,4/5,4/6,4/7,4/9,4/10,
  4/12,4/13,4/14,5/6,5/7,5/8,
  5/9,5/10,5/11,5/12,5/13,5/15,
  5/16,6/7,6/9,6/10,6/11,6/14,
  6/15,6/16,6/17,7/13,7/14,7/15,
  7/17,7/18,8/13,8/16,8/17,9/13,
  9/14,9/15,9/16,10/13,10/14,10/15,
  10/16,11/13,11/14,11/17,11/18,12/13,
  12/15,12/16,12/17,13/14,13/15,13/16,
  14/15,14/16,14/17,16/18,17/18}
  \draw[seedother] (v\a)--(v\b);
\foreach \a/\b in {2/3,2/4,2/15,2/16,3/4,3/15,
  3/17,4/16,4/17,7/8,7/9,7/10,
  7/11,8/9,8/10,8/12,9/11,9/12,
  10/11,10/12,11/12,15/16,15/17,16/17}
  \draw[seedremoved] (v\a)--(v\b);
\foreach \a/\b in {2/9,2/10,4/8,4/11,8/15,9/17,
  10/17,11/15}
  \draw[seedretained] (v\a)--(v\b);
\node[seedvertex,fill=orange!80!black] at (v1) {};
\node[font=\small,inner sep=0pt] at (90.000000:2.64) {$1$};
\node[seedvertex,fill=orange!80!black] at (v18) {};
\node[font=\small,inner sep=0pt] at (71.052632:2.64) {$18$};
\node[seedvertex] at (v7) {};
\node[font=\small,inner sep=0pt] at (52.105263:2.64) {$7$};
\node[seedvertex] at (v12) {};
\node[font=\small,inner sep=0pt] at (33.157895:2.64) {$12$};
\node[seedvertex] at (v9) {};
\node[font=\small,inner sep=0pt] at (14.210526:2.64) {$9$};
\node[seedvertex] at (v10) {};
\node[font=\small,inner sep=0pt] at (-4.736842:2.64) {$10$};
\node[seedvertex] at (v8) {};
\node[font=\small,inner sep=0pt] at (-23.684211:2.64) {$8$};
\node[seedvertex] at (v11) {};
\node[font=\small,inner sep=0pt] at (-42.631579:2.64) {$11$};
\node[seedvertex] at (v2) {};
\node[font=\small,inner sep=0pt] at (-61.578947:2.64) {$2$};
\node[seedvertex] at (v17) {};
\node[font=\small,inner sep=0pt] at (-80.526316:2.64) {$17$};
\node[seedvertex] at (v4) {};
\node[font=\small,inner sep=0pt] at (-99.473684:2.64) {$4$};
\node[seedvertex] at (v15) {};
\node[font=\small,inner sep=0pt] at (-118.421053:2.64) {$15$};
\node[seedvertex] at (v3) {};
\node[font=\small,inner sep=0pt] at (-137.368421:2.64) {$3$};
\node[seedvertex] at (v16) {};
\node[font=\small,inner sep=0pt] at (-156.315789:2.64) {$16$};
\node[seedvertex] at (v5) {};
\node[font=\small,inner sep=0pt] at (-175.263158:2.64) {$5$};
\node[seedvertex] at (v14) {};
\node[font=\small,inner sep=0pt] at (-194.210526:2.64) {$14$};
\node[seedvertex] at (v6) {};
\node[font=\small,inner sep=0pt] at (-213.157895:2.64) {$6$};
\node[seedvertex] at (v13) {};
\node[font=\small,inner sep=0pt] at (-232.105263:2.64) {$13$};
\node[seedvertex] at (v0) {};
\node[font=\small,inner sep=0pt] at (-251.052632:2.64) {$0$};
\draw[orange!80!black,line width=.45pt] (91.700000:2.86)--(91.700000:2.94)
  arc[start angle=91.700000,end angle=69.352632,radius=2.94cm]--(69.352632:2.86);
\node[font=\small,text=orange!80!black] at (80.526316:3.18) {$D$};
\draw[black!65,line width=.45pt] (53.805263:2.86)--(53.805263:2.94)
  arc[start angle=53.805263,end angle=31.457895,radius=2.94cm]--(31.457895:2.86);
\node[font=\small,text=black!65] at (42.631579:3.18) {$U$};
\draw[black!65,line width=.45pt] (15.910526:2.86)--(15.910526:2.94)
  arc[start angle=15.910526,end angle=-6.436842,radius=2.94cm]--(-6.436842:2.86);
\node[font=\small,text=black!65] at (4.736842:3.18) {$X$};
\draw[black!65,line width=.45pt] (-21.984211:2.86)--(-21.984211:2.94)
  arc[start angle=-21.984211,end angle=-44.331579,radius=2.94cm]--(-44.331579:2.86);
\node[font=\small,text=black!65] at (-33.157895:3.18) {$W$};
\draw[black!65,line width=.45pt] (-59.878947:2.86)--(-59.878947:2.94)
  arc[start angle=-59.878947,end angle=-82.226316,radius=2.94cm]--(-82.226316:2.86);
\node[font=\small,text=black!65] at (-71.052632:3.18) {$Y$};
\draw[black!65,line width=.45pt] (-97.773684:2.86)--(-97.773684:2.94)
  arc[start angle=-97.773684,end angle=-120.121053,radius=2.94cm]--(-120.121053:2.86);
\node[font=\small,text=black!65] at (-108.947368:3.18) {$Z$};
\draw[black!65,line width=.45pt] (-135.668421:2.86)--(-135.668421:2.94)
  arc[start angle=-135.668421,end angle=-158.015789,radius=2.94cm]--(-158.015789:2.86);
\node[font=\small,text=black!65] at (-146.842105:3.18) {$V$};
\end{tikzpicture}
\par\small (c) $r=14$, $n=19$
\end{minipage}
\par\medskip
\begin{tikzpicture}[font=\small]
\draw[red!75!black,dashed,line width=.65pt] (-6.6,0)--(-6,0);
\node[right] at (-5.9,0) {selected $R_r$};
\draw[blue!75!black,line width=.8pt] (-2.8,0)--(-2.2,0);
\node[right] at (-2.1,0) {selected $S_r$};
\draw[black!30,line width=.25pt] (1,0)--(1.6,0);
\node[right] at (1.7,0) {other edges};
\node[circle,draw=none,fill=orange!80!black,minimum size=1.2mm,inner sep=0pt] at (5,0) {};
\node[right] at (5.15,0) {deficient pair $D$};
\end{tikzpicture}
\caption{The three small modified bases used in
Lemma~\ref{lem:modified-group-base}, each connected and with magic constant $0$.
(a) over $\mathbb Z_{13}$, with degree $4$ at $D=\{6,7\}$ and degree $6$ elsewhere;
(b) over $\mathbb Z_{15}$, with degree $8$ at $D=\{3,12\}$ and degree $10$ elsewhere;
(c) over $\mathbb Z_{19}$, with degree $8$ at $D=\{1,18\}$ and degree $14$ elsewhere.
Outer brackets identify complementary pairs. Dashed red edges form $R_r$,
and solid blue edges form $S_r$; together they form the selected special copy of $H_r$.
All other edges are gray. The deficient pairs are marked in orange and are disjoint
from the selected copies. In (c), $U=\{7,12\}$, $V=\{3,16\}$,
$W=\{8,11\}$, $X=\{9,10\}$, $Y=\{2,17\}$, and $Z=\{4,15\}$.
The additional condition $N_B(x)\cap D=\varnothing$ holds for every $x\in X\cup Z$.
None of the selected copies is induced. These graphs are obtained using
the MILP models in the Appendix with the additional conditions specified there.}
\label{fig:small-modified-bases}
\end{figure}
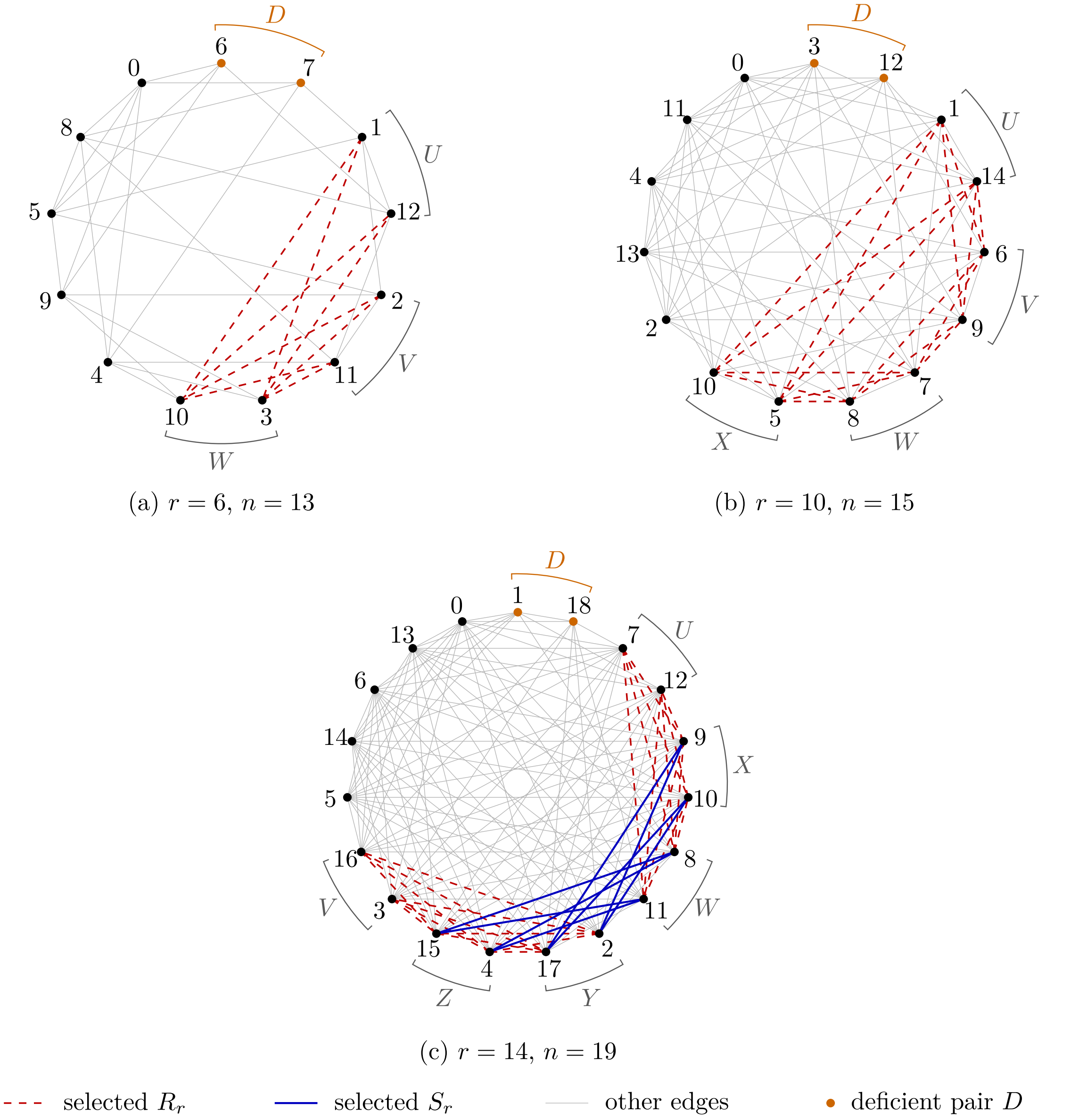
\begin{proof}[Proof of Theorem~\ref{thm:main}]
Retain the notation $\Gamma$, $K$, and $\Omega$ from Section~2.
By Lemma~\ref{lem:pairs}, the set $\Omega$ is partitioned into
\[
 M=\frac{n(|A|-1)}2
\]
two-element inverse pairs. In each case we start with a connected
graph of order $n$ labeled bijectively by $K$, with all weights zero.
We then use the inverse pairs in $\Omega$ to label the new vertices.

For $r=4$, take the graph in Corollary~\ref{lem:special-base} and
embed its labels into $K$ by $x\mapsto(0,x)$.
Apply Lemma~\ref{lem:step} once for each of the $M$ inverse pairs.
The graph remains connected and $4$-regular throughout.

Suppose next that $r\ge6$ and $|A|\equiv1\pmod4$.
Take the graph in Corollary~\ref{lem:special-base} and embed its
labels into $K$. The number $M$ is even, so we can apply
Lemma~\ref{lem:step} exactly $M/2$ times, using two unused inverse
pairs at each step. This includes $|A|=1$, for which $M=0$ and
no step is needed.

It remains to consider $r\ge6$ and $|A|\equiv3\pmod4$.
For $r=8,12$, use Corollary~\ref{lem:special-base} and embed the
labels into $K$. For $r=6,10,14$, take the connected graph and
$\Z_n$-distance magic labeling from Lemma~\ref{lem:modified-group-base}
and embed its labels into $K$. This includes the small bases
$(r,n)=(6,13),(10,15),(14,19)$ in
Figure~\ref{fig:small-modified-bases}. The specified copy of $H_r$
is special, and the deficient pair $D$ is complementary.

Set $s=1$ for $r=6,10,14$ and $s=3$ for $r=8,12$.
Since $n$ is odd and $|A|\equiv3\pmod4$, the number $M$ is odd.
Moreover, $M\ge n\ge s$ under the hypotheses of the theorem.
Reserve $s$ inverse pairs for the final step and apply
Lemma~\ref{lem:step} exactly $(M-s)/2$ times using the remaining
pairs. The selected copy of $H_r$ continues to consist of vertices
of degree $r$. In the cases with a deficient pair $D$, that pair
remains outside the selected copy and its neighbourhood is unchanged.

For $r=14$, Lemma~\ref{lem:modified-group-base} supplies a copy
with no edges from its pairs $X,Z$ to $D$. After each recursive step, the two
new pairs occupy the positions $X,Z$ in the selected replacement
copy, by Lemma~\ref{lem:step}. They have no neighbours in $D$
because the new pairs are joined only to each other and to the
previous copy. Thus the nonadjacency required in
Lemma~\ref{lem:final} holds before the final step, whether or not
any recursive step has been performed.
Apply that lemma using the $s$ reserved inverse pairs.

In every case, the initial graph is connected and each operation
preserves connectedness. The final graph is $r$-regular and all
its weights are zero. Its labels consist of every element of $K$
and every element of $\Omega$ exactly once. The labeling is
therefore a bijection onto $\Gamma$, and the graph has $n|A|$ vertices.
\end{proof}

\section{Concluding remarks}
The bounds in Theorem~\ref{thm:main} depend on the initial graphs
and the modifications used in the proof. The existence question
remains open for the cases with smaller odd values of $n$ not
covered by the theorem and for even degrees greater than $14$.
Extending the method to groups of even order would require
a separate treatment of any unused nonzero involutions, since
these elements cannot be partitioned into two-element inverse pairs.

The graphs constructed here need not admit an ordinary distance
magic labeling. For example, Figure~\ref{fig:Z15-distance-magic}
gives a connected $4$-regular $\mathbb Z_{15}$-distance magic graph
with magic constant $0$, whereas no $4$-regular distance magic
graph of order $15$ exists~\cite{KFK4reg}.

For $r\in\{4,8,10,12\}$, Theorem~\ref{thm:main} applies to every
finite Abelian group of odd order having a cyclic direct factor
of order at least $15$. By the structure theorem for finite
Abelian groups, the remaining nontrivial groups are precisely
those listed in the following problem.

\Needspace{16\baselineskip}
\begin{problem}
Let $r\in\{4,8,10,12\}$. Determine which finite Abelian groups
\[
\Gamma\cong
(\mathbb Z_3)^{m_1}\times(\mathbb Z_9)^{m_2}
\quad\text{or}\quad
\Gamma\cong(\mathbb Z_p)^m,
\]
where $m_1,m_2\ge0$, $m_1+m_2\ge1$, $m\ge1$, and
$p\in\{5,7,11,13\}$, admit a connected $r$-regular
$\Gamma$-distance magic graph with magic constant $0$.
\end{problem}

\textbf{Generative AI statement.} The authors used ChatGPT (OpenAI) for language editing and assistance with the preparation of figures. All proofs were developed by the authors. The authors reviewed and edited the AI-assisted material and take full responsibility for the content of the manuscript.
\FloatBarrier

\section*{Appendix}

\subsection*{MILP formulation}

Let $n,r,k\in\mathbb{N}_+$ and let
$V=\{0,1,\ldots,n-1\}$. Our goal is to find a simple
$r$-regular graph $G=(V,E)$ such that
\begin{equation}
    \sum_{v\in N_G(u)}v\equiv 0\pmod n
    \qquad \text{for every }u\in V.
    \label{eq:milp_modular_condition}
\end{equation}
Additionally, we require $k$ pairwise disjoint subsets
$U_1,\ldots,U_k\subseteq V$, each of which is a
complementary pair, i.e., $U_i=\{u,v\}$ for some
distinct vertices $u,v\in V$ satisfying $u+v\equiv0\pmod n$.
Let
$\mathcal{A}\subseteq
    \bigl\{\{i,j\}:1\leq i<j\leq k\bigr\}$
be a fixed family specifying the required joins $U_i U_j$.
For every $\{i,j\}\in\mathcal{A}$, we require that
every vertex of $U_i$ is adjacent to every vertex
of $U_j$.

We formulate this problem as a Mixed Integer Linear
Programming (MILP) feasibility problem.

We introduce the following variables:
\begin{itemize}
    \item $e_{uv}\in\{0,1\}$ for $u,v\in V$,
    indicating whether $\{u,v\}$ is an edge of $G$;
    \item $I_{iv}\in\{0,1\}$ for $i\in\{1,\ldots,k\}$
    and $v\in V$, indicating whether $v\in U_i$;
    \item $q_u\in\mathbb{Z}$ for $u\in V$,
    used to express the modular condition.
\end{itemize}

The MILP model is given by
\begin{align}
    \min\quad & 0 \notag\\
    \text{subject to}\quad
    & e_{uu}=0
    && \forall u\in V,
    \label{eq:milp_no_loops}\\
    & e_{uv}=e_{vu}
    && \forall u,v\in V,\ u<v,
    \label{eq:milp_symmetry}\\
    & \sum_{v\in V - \{u\}}e_{uv}=r
    && \forall u\in V,
    \label{eq:milp_degree}\\
    & \sum_{v\in V}v e_{uv}=nq_u
    && \forall u\in V,
    \label{eq:milp_modular}\\
    & \sum_{i=1}^{k}I_{iv}\leq 1
    && \forall v\in V,
    \label{eq:milp_disjoint}\\
    & \sum_{v\in V}I_{iv}=2
    && \forall i\in\{1,\ldots,k\},
    \label{eq:milp_pair_size}\\
    & \sum_{v\in V}v I_{iv}=n
    && \forall i\in\{1,\ldots,k\},
    \label{eq:milp_complementary}\\
    & e_{uv}\geq I_{iu}+I_{jv}-1
    && \substack{
       \forall\{i,j\}\in\mathcal{A},\ i<j,\\
       \forall (u,v)\in V^2,\ u\neq v,
       }
    \label{eq:milp_joins}
\end{align}

Constraints~\eqref{eq:milp_no_loops} and
\eqref{eq:milp_symmetry} ensure that $G$ is a simple
undirected graph, while \eqref{eq:milp_degree} enforces
$r$-regularity. Constraint~\eqref{eq:milp_modular}
guarantees that the sum of the labels of the neighbors
of each vertex is divisible by $n$. Constraints~\eqref{eq:milp_disjoint} and
\eqref{eq:milp_pair_size} ensure that the sets
$U_1,\ldots,U_k$ are pairwise disjoint and each
contains exactly two vertices. Moreover,
\eqref{eq:milp_complementary} guarantees that the
two vertices in each set form a complementary pair.
Finally, \eqref{eq:milp_joins} forces all edges
between $U_i$ and $U_j$ whenever
$\{i,j\}\in\mathcal{A}$.

Thus, a graph satisfying the above requirements
exists if and only if the MILP model is feasible.
If a graph obtained from this model is connected,
then it is a proper base graph.

\paragraph{Modification for two vertices of a different degree.}

Suppose that we want to find a graph with exactly
two vertices of degree $d$ and all remaining vertices
of degree $r$, where $d\neq r$.

For each $u\in V$, we introduce two additional
binary variables $\alpha_u$ and $\beta_u$.
The variable $\alpha_u$ indicates whether $u$
has degree $d$, while $\beta_u$ indicates whether
$u$ has degree $r$.

We replace constraint~\eqref{eq:milp_degree} with
the following constraints:
\begin{align*}
    \alpha_u+\beta_u&=1
    && \forall u\in V,\\
    \sum_{u\in V}\alpha_u&=2,
    \\
    \sum_{v\in V - \{u\}}e_{uv}
    &=d\alpha_u+r\beta_u
    && \forall u\in V.
\end{align*}

If we want to ensure that no vertex of degree $d$
belongs to a given set $U_i$, it is sufficient
to add
\begin{equation*}
    I_{iv}+\alpha_v\leq 1
    \qquad \forall v\in V.
\end{equation*}

If we want to forbid all edges between
vertices of a given set $U_i$ and vertices of
degree $d$, we can add
\begin{equation*}
    e_{uv}+I_{iu}+\alpha_v\leq 2
    \qquad \forall (u,v)\in V^2,\ u\neq v.
\end{equation*}

For the three modified bases in Figure~\ref{fig:small-modified-bases},
take
\[
 (n,r,d,k)\in\{(13,6,4,3),(15,10,8,4),(19,14,8,6)\}
\]
and let $\mathcal A$ specify all joins of the corresponding $H_r$
in Table~\ref{tab:patterns}. Impose the exclusion constraint
$I_{iv}+\alpha_v\leq1$ for every prescribed pair and every vertex,
so that $D$ is disjoint from the selected copy. Also impose
\[
 \sum_{v\in V}v\alpha_v=n.
\]
Since exactly two vertices have $\alpha_v=1$, this equality makes
$D$ complementary. For $r=14$, impose the nonadjacency constraint
above for the pairs in the positions $X$ and $Z$.
The displayed graphs have all the properties asserted in
Lemma~\ref{lem:modified-group-base}, including connectedness,
as verified in its proof.


\begin{thebibliography}{9}
\bibitem{ArumugamSurvey}
S.~Arumugam, D.~Fron\v{c}ek and N.~Kamatchi,
\emph{Distance Magic Graphs---A Survey},
J.~Indones.\ Math.\ Soc., Special Edition (2011), 11--26.
\bibitem{FroncekGroup}
D.~Fron\v{c}ek, \emph{Group distance magic labeling of Cartesian product of cycles},
Australas.\ J.\ Combin.\ \textbf{55} (2013), 167--174.

\bibitem{FroncekFair}
D.~Fron\v{c}ek, P.~Kov\'a\v{r} and T.~Kov\'a\v{r}ov\'a,
\emph{Fair incomplete tournaments},
Bull.\ Inst.\ Combin.\ Appl.\ \textbf{48} (2006), 31--33.
\bibitem{Gallian}
J.~A.~Gallian, \emph{A dynamic survey of graph labeling},
Electron.\ J.\ Combin.\ DS6 (2022).
\bibitem{KFK4reg}
P.~Kov\'a\v{r}, D.~Fron\v{c}ek and T.~Kov\'a\v{r}ov\'a,
\emph{A note on $4$-regular distance magic graphs},
Australas.\ J.\ Combin.\ \textbf{54} (2012), 127--132.
\bibitem{KovKrb14}
P.~Kov\'a\v{r} and M.~Krbe\v{c}ek,
\emph{On $14$-regular distance magic graphs},
Electron.\ J.\ Graph Theory Appl.\ \textbf{12}(1) (2024), 35--41.

\bibitem{KSKK6810}
P.~Kov\'a\v{r}, A.~Silber, P.~Kabel\'{\i}kov\'a--Hru\v{s}kov\'a and M.~Krav\v{c}enko,
\emph{On regular distance magic graphs of odd order},
J.\ Combin.\ Math.\ Combin.\ Comput.\ \textbf{117} (2023), 55--64.

\end{thebibliography}
\end{document}